\documentclass[11pt]{article}
\usepackage{amsfonts}
\usepackage{amssymb}
\usepackage{amsthm}
\usepackage{amsmath}
\usepackage{graphicx}
\usepackage{empheq}
\usepackage{indentfirst}
\usepackage{cite} 
\usepackage{mathrsfs}
\usepackage{cases}
\usepackage{graphics}
\usepackage{xcolor}  
\usepackage{amsmath,bm}
\newtheorem{theorem}{\textbf{Theorem}}[section]
\newtheorem{lemma}{\textbf{Lemma}}[section]
\newtheorem{proposition}{\textbf{Proposition}}[section]
\newtheorem{corollary}{\textbf{Corollary}}[section]
\newtheorem{remark}{\textbf{Remark}}[section]
\newtheorem{definition}{\textbf{Definition}}[section]
\allowdisplaybreaks[4] 
\def\be{\begin{equation}}
\def\ee{\end{equation}}
\def\bea{\begin{eqnarray}}
\def\eea{\end{eqnarray}}
\def\bt{\begin{theorem}}
\def\et{\end{theorem}}
\def\bl{\begin{lemma}}
\def\el{\end{lemma}}
\def\br{\begin{remark}}
\def\er{\end{remark}}
\def\bp{\begin{proposition}}
\def\ep{\end{proposition}}
\def\bc{\begin{corollary}}
\def\ec{\end{corollary}}
\def\bd{\begin{definition}}
\def\ed{\end{definition}}
\begin{document}

\title{Global Weak Solutions to a Cahn--Hilliard--Navier--Stokes System with Chemotaxis and Singular Potential}

\author{
{Jingning He}\footnote{School of Mathematical Sciences, Fudan University, Handan Road 220, Shanghai 200433, China. Email:  \textit{jingninghe2020@gmail.com}
}}

\date{\today}

\maketitle


\begin{abstract}
\noindent We analyze a diffuse interface model that describes the dynamics of incompressible two-phase flows with chemotaxis effect.
The PDE system couples a Navier--Stokes equation for the fluid velocity, a convective Cahn--Hilliard equation for the phase field variable with an advection-diffusion-reaction equation for the nutrient density. For the system with a singular potential, we prove the existence of global weak solutions in both two and three dimensions. Besides, in the two dimensional case, we establish a continuous dependence result that implies the uniqueness of global weak solutions. The singular potential guarantees that the phase field variable always stays in the physically relevant interval $[-1,1]$ during time evolution. This property enables us to obtain the well-posedness result without any extra assumption on the coefficients that has been made in the previous literature. \medskip 
\end{abstract}

\section{Introduction}
\setcounter{equation}{0}
\noindent 
Diffuse interface models have emerged as an efficient mathematical tool describing the complex dynamics of mixtures in materials science \cite{CH}, fluid dynamics \cite{Gur,HH,J,LS,LT98}, and mathematical biology, e.g., the tumor growth process \cite{CLLW,GLSS,HvO,OHP}. In the diffuse interface framework, large interface deformations and topological changes of the interfaces of the mixture can be handled naturally. 

In this paper, we consider a Cahn--Hilliard--Navier--Stokes type system
\begin{subequations}
\begin{alignat}{3}
&\partial_t  \bm{ v}+(\bm{ v} \cdot \nabla)  \bm {v}-\textrm{div} (  2\eta(\varphi) D\bm{v} )+\nabla p=(\mu+\chi \sigma)\nabla \varphi,\qquad &\textrm{in}\ \Omega\times(0,T),\label{f3.c} \\
&\textrm{div}\ \bm{v}=0,\ &\textrm{in}\ \Omega\times(0,T),\label{f3.c1}\\
&\partial_t \varphi+(\bm{v} \cdot \nabla) \varphi=\Delta \mu-\alpha(\varphi-c_0),\ &\textrm{in}\ \Omega\times(0,T),\label{f1.a} \\
&\mu=A\varPsi'(\varphi)-B\Delta \varphi-\chi \sigma, \ &\textrm{in}\ \Omega\times(0,T),\label{f4.d} \\
&\partial_t \sigma+(\bm{v} \cdot \nabla) \sigma= \Delta (\sigma+\chi(1-\varphi))-\mathcal{C} h(\varphi) \sigma +S,\ &\textrm{in}\ \Omega\times(0,T), \label{f2.b}  
\end{alignat}
\end{subequations}
subject to the following boundary conditions
\begin{alignat}{3}
&\bm{v}=\mathbf{0},\quad {\partial}_{\bm{n}}\varphi={\partial}_{\bm{n}}\sigma={\partial}_{\bm{n}}\mu=0,\qquad\qquad &\textrm{on}& \   \partial\Omega\times(0,T),
\label{boundary}
\end{alignat}
as well as initial conditions
\begin{alignat}{3}
&\bm{v}(0)=\bm{v}_{0},\ \ \varphi(0)=\varphi_{0},\ \ \sigma(0)=\sigma_{0}, \qquad &\textrm{in}&\ \Omega.
\label{ini0}
\end{alignat}
Here, $\Omega \subset\mathbb{R}^d\ (d=2,\ 3)$ is a bounded domain with smooth boundary $\partial\Omega$ and $\bm{n}=\bm{n}(x)$ denotes the unit outward normal vector on $\partial\Omega$. $S=S(x,t)$ is a given function standing for possible external source.

The system \eqref{f3.c}--\eqref{f2.b} can be viewed as a simplified version of the general thermodynamically consistent diffuse interface model that was derived in Lam and Wu \cite{LW} for a two-component incompressible fluid mixture with a chemical species subject to diffusion as well as other transport mechanisms like convection and chemotaxis (see the earlier work \cite{Sitka} for a more specific setting in the context of tumor growth modelling such that the mixture describes a tumor surrounded by healthy tissues). The order parameter $\varphi$ denotes the difference in volume fractions of the mixture such that the region $\{\varphi=1\}$ represents fluid 1 and $\{\varphi=-1\}$ represents fluid 2 (i.e., the values $\pm 1$ represent the pure concentrations). The fluid velocity $\bm{v}$ is taken as the volume-averaged velocity with $D\bm{v}=\frac{1}{2}(\nabla\bm{ v}+(\nabla\bm{ v}) ^ \mathrm{T})$ being the symmetrized velocity gradient, and the scalar function $p$ is the (modified) pressure. The variable $\sigma$ denotes the concentration of the chemical species (e.g., nutrient) and $\mu$ stands for the chemical potential associated to $(\varphi, \sigma)$. Equations \eqref{f3.c} and \eqref{f3.c1} represent the momentum balance for the incompressible fluid mixture, while equations \eqref{f1.a} and \eqref{f4.d} constitute a convective Cahn--Hilliard system for the order parameter $\varphi$, and equation \eqref{f2.b} is an advection-diffusion-reaction equation for the chemical density $\sigma$.

For the sake of simplicity, in this paper we assume that the density difference of the mixture as well as the mass transfer between the two components are negligible. Besides, we assume that the mobilities are positive constants (set to be $1$). The source term in the Cahn--Hilliard equation \eqref{f1.a} may correspond to some biological mechanisms like proliferation, apoptosis of cells in the tumor growth modelling. Here, we only take a simple form, i.e., of Oono's type $-\alpha(\varphi-c_0)$ (cf. \cite{GGM,Mi11}), where $\alpha\geq 0$, $c_0\in(-1,1)$. For further discussions on biologically relevant mass source terms of the Cahn--Hilliard equation, we refer to \cite{Fa15,Fa17}. The nutrient consumption is prescribed by the term $\mathcal{C}h(\varphi)\sigma$, where the non-negative constant $\mathcal{C}$ represents the consumption rate and the function $h$ is an interpolation with $h(-1)=0$ and $h(1)=1$, for instance, the simplest choice could be $h(\varphi) = \frac{1}{2}(1+\varphi)$ (cf. \cite{GLSS}). We allow that the binary fluid mixture may have unmatched viscosities. Assuming that $\eta_1$, $\eta_2>0$ are viscosities of the two homogeneous fluids, viscosity of the mixture can be modeled by the concentration dependent term $\eta=\eta(\varphi)$, for instance, a typical form is the linear combination (see, e.g., \cite{J}):
\be
\eta(r)=\eta_1\frac{1+r}{2}+\eta_2\frac{1-r}{2}, \quad \forall\, r\in[-1,1].
\label{vis}
\ee  

In \eqref{f4.d}, the positive constants $A, B$ are related to the surface tension and the thickness of the interfacial layers (i.e., the diffuse interface). The nonlinear function $\varPsi'$ denotes the derivative of a potential $\varPsi$ that has a double-well structure, with two minima and a local unstable maximum in between. A physically significant example is the following logarithmic type:
 \be
\varPsi (r)=\frac{\theta}{2}[(1-r)\ln(1-r)+(1+r)\ln(1+r)]+\frac{\theta_{c}}{2}(1-r^2),\quad \forall\, r\in(-1,1),
\label{pot}
 \ee
 with $0<\theta<\theta_{c}$ (see e.g., \cite{CH,CMZ}). It is referred to as a singular potential since its derivative $\varPsi'$ blows up at the pure phases $\pm 1$. In the literature, the singular potential $\varPsi$ is often approximated by a fourth-order polynomial
  \be
\varPsi(r)=\frac{1}{4}(1-r^2)^2,\quad r\in\mathbb{R}, \label{regular}
 \ee
or some more general polynomial function.
 
The coupling structure of system \eqref{f3.c}--\eqref{f2.b} is reflected in terms of the capillary force $(\mu+\chi \sigma)\nabla \varphi$ (indeed only depending on $\varphi$ in view of \eqref{f4.d}), the viscous stress tensor with a concentration dependent viscosity $\eta(\varphi)$,  
the advection terms $\bm{v}\cdot \nabla \varphi$, $\bm{v}\cdot \nabla \sigma$, and two extra terms involving the parameter $\chi$. In particular, the constant coefficient $\chi$ is related to certain specific transport mechanisms such as chemotaxis and active transport in the context of tumor growth modelling (see e.g., \cite{GLSS, GL17e}). To see this, we reformulate equations \eqref{f1.a} and \eqref{f2.b} as 
\begin{align*}
&\partial_t \varphi+(\bm{v} \cdot \nabla) \varphi+\mathrm{div}\bm{q}_\varphi=0, \qquad \partial_t \sigma+(\bm{v} \cdot \nabla) \sigma+ \mathrm{div}\bm{q}_\sigma=-\mathcal{C} h(\varphi) \sigma +S,\\
&\text{with fluxes}\ \ \bm{q}_\varphi:= -\nabla \mu= -\nabla(A\varPsi'(\varphi)-B\Delta \varphi-\chi \sigma),\quad
\bm{q}_\sigma:=-\nabla(\sigma- \chi\varphi).
\end{align*}
For $\chi\geq 0$, the term $\chi \nabla \sigma$ in $\bm{q}_\varphi$ represents the chemotactic response to the nutrient (i.e., movement of fluid 2 towards regions of high $\sigma$), while the other term $\chi\nabla \varphi$ in $\bm{q}_\sigma$ drives the chemical
species towards fluid 2, i.e., the region $\{\varphi = 1\}$, leading to a persistent concentration difference between the mixture components against the usual diffusion effect (especially near the interface where $\nabla \varphi$ is non-zero). In spite of the complicated coupling structure, we observe that the system \eqref{f3.c}--\eqref{f2.b} admits a basic energy law 
\begin{align}
& \frac{d}{dt} \int_{\Omega} \Big[ \frac{1}{2}|\bm{v}|^2+A\varPsi(\varphi)+\frac{B}{2}|\nabla \varphi|^2+\frac{1}{2}|\sigma|^2+\chi\sigma(1-\varphi) \Big] dx \notag \\
& \qquad +\int_{\Omega} \Big[ 2\eta(\varphi)|D\bm{v}|^2 + |\nabla \mu|^2+|\nabla(\sigma+\chi(1-\varphi))|^2\Big] dx\nonumber\\
&\quad  =\int_\Omega \left[-\alpha(\varphi-c_0)\mu+(-\mathcal{C} h(\varphi) \sigma +S)(\sigma +\chi(1-\varphi))\right] dx, 
\label{BEL}
\end{align}
which plays an important role in the study of its global well-posedness.

To the best of our knowledge, the only known analytic result for problem \eqref{f3.c}--\eqref{ini0} is due to Lam and Wu \cite{LW} (with described source terms).
Under the choice of a regular potential including the prototype \eqref{regular}, 
they establish the existence of global weak solutions in two and three dimensions for prescribed mass transfer terms as well as the existence and uniqueness of global strong solutions in two dimensions. However, it is worth noting that the fourth-order Cahn--Hilliard equation for $\varphi$ does not have a maximum principle, with a regular potential like \eqref{regular} one cannot guarantee the solution $\varphi$ to take values in the physical interval $[-1, 1]$ (see e.g., \cite[Remark 2.1]{CMZ} for a counterexample). Due to this  technical issue, in order to prove the existence of global weak solutions, the authors of \cite{LW} have to impose the following assumption on the coefficients $A$ and $\chi$:
\begin{align}
A>\frac{2\chi^2}{C_3},
\label{asscoe}
\end{align}
where $C_3$, $C_4$ are positive constants such that $\Psi(r)\geq C_3r^2-C_4$ for $r\in \mathbb{R}$ (see also \cite{GL17e} for a similar situation in the fluid-free case with more general mass source terms). The assumption \eqref{asscoe} arises from using H\"{o}lder's and Young's inequalities to control the crossing term $\int_\Omega \chi\sigma(1-\varphi) dx$ in the free energy (see \eqref{BEL}).  
Besides, the uniqueness of weak solution to problem \eqref{f3.c}--\eqref{ini0} in the two dimensional case was unsolved in \cite{LW}.  

Our aim in this paper is to expand the recent analysis for the initial boundary value problem \eqref{f3.c}--\eqref{ini0}. Taking  a singular potential into account (e.g., the physically relevant logarithmic type \eqref{pot}) and without using the above restrictive assumption \eqref{asscoe} on coefficients, we are able to prove:
\begin{itemize}
\item[(1)] existence of global weak solutions to problem \eqref{f3.c}--\eqref{ini0} in both two and three dimensions (see Theorem \ref{main});
\item[(2)] a continuous dependence result in dimension two (see Theorem \ref{the2}) that also yields the uniqueness of global weak solution (see Corollary \ref{uniq}).
\end{itemize}

The assumption \eqref{asscoe} seems to imply that in order to obtain the global weak solution etc of problem \eqref{f3.c}--\eqref{ini0}, the effects of chemotaxis as well as active transport cannot be too strong. From the technical point of view, without this assumption, it is not clear whether the Galerkin approximation scheme used in \cite{LW} still works for singular potentials (see hypotheses (H2) in the next section). This is because after a regularization of the singular potential $\varPsi$ (see e.g., \cite{B,MT}), the approximate ansatz for $\varphi$ in the usual Galerkin scheme does not belong to the interval $[-1,1]$ due to the lack of maximum principle, which yields difficulties to derive necessary uniform \textit{a priori} estimates. To overcome this difficulty, we shall make use of an alternative method with a semi-Galerkin scheme, that is, performing a Galerkin approximation only for the Navier--Stokes equation of $\bm{v}$, but keeping the equations for the other variables $(\varphi, \sigma)$ and then applying a fixed point argument. Since the Cahn--Hilliard equation with singular potential is solved separately in this procedure, taking advantage of the existing literature (cf. e.g., \cite{A2009,A2007,MT}), we can guarantee the property $\varphi \in [-1,1]$ for approximate solutions. This approach has been successfully applied to other nonlinear coupled systems, for instance, the Ericksen--Leslie system for incompressible liquid crystal flow \cite{LfL} and a diffuse interface model for incompressible binary fluids with thermal Marangoni effect \cite{W2017}. We remark that the property $\varphi \in [-1,1]$ is also important in view of the variable viscosity $\eta(\varphi)$, defined via the relation \eqref{vis}, which can become negative if $\varphi$ is outside the physical interval $[-1, 1]$. On the other hand, the non-constant viscosity, non conservation of mass  and the coupling with nutrient equation lead to additional mathematical difficulties to prove the uniqueness of global weak solutions of problem \eqref{f3.c}--\eqref{ini0} in two dimensions. We shall extend the method introduced in the recent work \cite{GMT} for the Cahn--Hilliard--Navier--Stokes system to derive a continuous dependence estimate with respect to initial data in weaker norms of the solution, from which the uniqueness of weak solutions follows. 

Before ending the introduction, let us give, without any claim of completeness, a brief overview of related mathematical analysis results in the literature. When the nutrient interaction is neglected and $\alpha=0$, system \eqref{f3.c}--\eqref{f4.d} reduces to the well-known Model H for the motion of incompressible, viscous two-phase flow \cite{HH, Gur}. The resulting Cahn--Hilliard--Navier--Stokes system with regular potentials has been widely studied, see for instance, \cite{BGM,GG2010,GG2010b,ZWH} and the references cited therein. For the Cahn--Hilliard--Navier--Stokes system with unmatched viscosities and the logarithmic potential, we refer to \cite{A2009,B,GMT}, see also \cite{GGM,GGW19} for the case with more involved boundary conditions accounting for the moving contact line, \cite{MT} for the Cahn--Hilliard--Oono--Navier--Stokes system (i.e., $\alpha>0$) with constant viscosity, and  \cite{A2012,A2013,GGW19} for fluid mixtures with different densities. On the other hand, when the fluid interaction in system \eqref{f3.c}--\eqref{f4.d} is neglected, we refer to \cite{GL17e} for the well-posedness of a Cahn--Hilliard type system with chemotaxis and active transport (see \cite{GL17} for the case with Dirichlet boundary conditions) and to \cite{MRS} for a first study on the long-time behavior without transport mechanisms. It is worth mentioning that diffuse interface models with other types of fluid interaction have also been extensively investigated in the literature, for instance, we refer to \cite{DFRSS,GLSS,GGW,HWW,JWZ,LTZ,WW,WZ} for the Cahn--Hilliard--Darcy system and to \cite{BCG,CG,EG19jde,EG19sima} for the Cahn--Hilliard--Brinkman system with various extensions especially in the recent study of tumor growth modelling.  

The remaining part of this paper is organized as follows. In Section \ref{pm}, we introduce the functional settings and state the main results. In Section \ref{ws}, we prove the existence of global weak solutions in both two and three dimensions. In Section \ref{uw}, we derive a continuous dependence result and prove the uniqueness of weak solution in dimension two. In the Appendix, we provide some details of the semi-Galerkin approximate scheme that is used in the proof of the existence result.


\section{Main Results}\label{pm}
\setcounter{equation}{0}
\subsection{Preliminaries}
Throughout the paper, we assume that $\Omega \subset\mathbb{R}^d$ ($d=2,3$) is a bounded domain with smooth boundary $\partial\Omega$. For the standard Lebesgue and Sobolev spaces, we use the notations $L^{p} := L^{p}(\Omega)$ and $W^{k,p} := W^{k,p}(\Omega)$ for any $p \in [1,+\infty]$, $k > 0$ equipped with the norms $\|\cdot\|_{L^{p}}$ and $\|\cdot\|_{W^{k,p}}$.  In the case $p = 2$ we use $H^{k} := W^{k,2}$ and the norm $\|\cdot\|_{H^{k}}$. The norm and inner product on $L^{2}(\Omega)$ are simply denoted by $\|\cdot\|$ and $(\cdot,\cdot)$, respectively. 
The dual space of a Banach space $X$ is denoted by $X'$, and the duality pairing between $X$ and its dual will be denoted by
$\langle \cdot,\cdot\rangle_{X}$. Given an interval $J$ of $\mathbb{R}^+$, we introduce
the function space $L^p(J;X)$ with $p\in [1,+\infty]$, which
consists of Bochner measurable $p$-integrable
functions with values in the Banach space $X$. 
The boldface letter $\bm{X}$ denotes the vectorial space $X^d$ endowed with the product structure.

For every $f\in H^1(\Omega)'$, we denote by $\overline{f}$ its  generalized mean value over $\Omega$ such that
$\overline{f}=|\Omega|^{-1}\langle f,1\rangle_{H^1}$; if $f\in L^1(\Omega)$, then its mean is simply given by $\overline{f}=|\Omega|^{-1}\int_\Omega f \,dx$. 
As the pressure function in \eqref{f3.c} is determined up to a time-dependent constant, we introduce the space $L^2_{0}(\Omega):=\{f\in L^2(\Omega):\overline{f} =0\}$. Besides, in view of the homogeneous Neumann boundary condition \eqref{boundary}, we also set
$H^2_{N}(\Omega):=\{f\in H^2(\Omega):\,\partial_{\bm{n}}f=0 \ \textrm{on}\  \partial \Omega\}$.  We will use the Poincar\'{e}--Wirtinger inequality \cite[Section 5.8.1]{E}:
\begin{equation}
\label{poincare}
\|f-\overline{f}\|\leq C_P\|\nabla f\|,\quad \forall\,
f\in H^1(\Omega),
\end{equation}
where $C_P$ is a constant depending only on $d$ and $\Omega$.
Consider the realization of the minus Laplacian with homogeneous Neumann boundary condition 
$\mathcal{A}_N\in \mathcal{L}(H^1(\Omega),H^1(\Omega)')$ defined by 
\begin{equation}\nonumber
   \langle \mathcal{A}_N u,v\rangle_{H^1} := \int_\Omega \nabla u\cdot \nabla v \, dx,\quad \text{for }\,u,v\in H^1(\Omega).
\end{equation}
Then for the linear spaces
$$
V_0=\{ u \in H^1(\Omega):\ \overline{u}=0\}, \quad
V_0'= \{ u \in H^1(\Omega)':\ \overline{u}=0 \},
$$
the restriction of $\mathcal{A}_N$ from $V_0$ onto $V_0'$
is an isomorphism. In particular, $\mathcal{A}_N$ is positively defined on $V_0$ and self-adjoint. We denote its inverse map by $\mathcal{N} =\mathcal{A}_N^{-1}: V_0'
\to V_0$. Note that for every $f\in V_0'$, $u= \mathcal{N} f \in V_0$ is the unique weak solution of the Neumann problem
$$
\begin{cases}
-\Delta u=f, \quad \text{in} \ \Omega,\\
\partial_{\bm{n}} u=0, \quad \ \  \text{on}\ \partial \Omega.
\end{cases}
$$
Besides, we have
\begin{align}
&\langle \mathcal{A}_N u, \mathcal{N} g\rangle_{V_0} =\langle  g,u\rangle_{V}, \quad \forall\, u\in V, \ \forall\, g\in V_0',\label{propN1}\\
&\langle  g, \mathcal{N} f\rangle_{V_0} 
=\langle f, \mathcal{N} g\rangle_{V_0} = \int_{\Omega} \nabla(\mathcal{N} g)
\cdot \nabla (\mathcal{N} f) \, dx, \quad \forall \, g,f \in V_0',\label{propN2}
\end{align}
and the chain rule
\begin{align}
&\langle \partial_t u, \mathcal{N} u(t)\rangle_{V_0}=\frac{1}{2}\frac{d}{dt}\|\nabla \mathcal{N} u\|^2,\ \ \textrm{a.e. in}\ (0,T),\nonumber
\end{align}
for any $u\in H^1(0,T; V_0')$. For any $f\in V_0'$, we set $\|f\|_{V_0'}=\|\nabla \mathcal{N} f\|$.
It is well-known that $f \to \|f\|_{V_0'}$ and $
f \to(\|f-\overline{f}\|_{V_0'}^2+|\overline{f}|^2)^\frac12$ are
equivalent norms on $V_0'$ and $H^1(\Omega)'$,
respectively. Besides, according to Poincar\'{e}'s inequality \eqref{poincare}, we
see that $f\to \|\nabla f\|$,  $f\to (\|\nabla f\|^2+|\overline{f}|^2)^\frac12$ are equivalent norms on $V_0$ and $H^1(\Omega)$.
We also report the following standard Hilbert interpolation inequality and elliptic estimates for the Neumann problem
\begin{align}
\|f\| &\leq \|f\|_{V_0'}^{\frac12} \| \nabla f\|^{\frac12},
\qquad \forall\, f \in V_0,\label{I}\\
\|\nabla \mathcal{N} f\|_{\bm{H}^{k}(\Omega)}& \leq C \|f\|_{H^{k-1}(\Omega)},
\qquad \forall\, f\in H^{k-1}(\Omega)\cap L^2_0(\Omega),\quad k\in\mathbb{N}.\label{N}
\end{align}
We also consider the operator $\mathcal{A}_1 := I-\Delta $ with homogeneous Neumann boundary condition that is an unbounded operator $L^2(\Omega)$ with domain $D(\mathcal{A}_1) =H^2_N(\Omega)$. It is well-known that
$\mathcal{A}_1$ is a positive, unbounded, self-adjoint operator in $L^2(\Omega)$ with a compact inverse (denoted by $\mathcal{N}_1:=\mathcal{A}^{-1}_1$), see, e.g., \cite[Chapter II, Section 2.2]{T}. Then $f \to \|\mathcal{N}_1^\frac12 f\|$ is also an equivalent norm on $H^1(\Omega)'$.
 
Next, we introduce the classical function spaces for the Navier--Stokes equations (see e.g., \cite{G,S}). For a vector-valued/tensor-valued Banach space $\bm{X}$, we denote $\bm{X}_{\mathrm{div}} $, $\bm{X}_{0,\mathrm{div}} $ by the closure of $C_{\mathrm{div}}^{\infty}(\Omega)=\{\bm{f}\in (C^{\infty}(\Omega))^{d}:\ \textrm{div}\bm{f}=0\}$,  $C_{0,\mathrm{div}}^{\infty}(\Omega)=\{\bm{f}\in (C_0^{\infty}(\Omega))^{d}:\ \textrm{div}\bm{f}=0\}$ with respect to the $\bm{X}$-norm, respectively. For $\bm{X}=\bm{L}^2(\Omega)$, we have the notation $\bm{L}^2_{0,\mathrm{div}}(\Omega)=\bm{L}^2_{\mathrm{div}}(\Omega)$. The space $\bm{H}^1_{0,\mathrm{div}}(\Omega)$ is equipped with the scalar product
 \be
(\bm{u},\bm{v})_{\bm{H}^1_{0,\mathrm{div}}}:=(\nabla \bm{u},\nabla \bm{v}),\quad   \forall\, \bm{u},\, \bm{v}  \in {\bm{H}^1_{0,\mathrm{div}}(\Omega)}.\nonumber
 \ee
 It is well known that $\bm{L}^2(\Omega)$ can be decomposed into $\bm{L}^2_{\mathrm{div}}(\Omega)\oplus\bm{G}(\Omega)$, where $\bm{G}(\Omega):=\{\bm{f}\in\bm{L}^2(\Omega): \exists\, z\in H^1(\Omega),\ \bm{f}=\nabla z\}$. Then for any function $\bm{f} \in \bm{L}^2(\Omega)$, there holds the Helmholtz--Weyl decomposition (see  \cite[Chapter \uppercase\expandafter{\romannumeral3}]{G}):
\be
\bm{f}=\bm{f}_{0}+\nabla z,\quad\text{where}\  \bm{f}_{0} \in \bm{L}^2_{\mathrm{div}}(\Omega),\ \nabla z \in \bm{G}(\Omega).\nonumber
 \ee
Consequently, we can define the Helmholtz--Leray projection onto the space of divergence-free functions $\bm{P}:\bm{L}^2(\Omega)\to \bm{L}^2_{\mathrm{div}}(\Omega)$ such that $\bm{P}(\bm{f})=\bm{f}_{0}$.  

We now invoke the Stokes operator $\bm{S}: \bm{H}^1_{0,\mathrm{div}}(\Omega)\cap\bm{H}^2(\Omega)\to\bm{L}^2_{\mathrm{div}}(\Omega)$ such that  
\be
(\bm{S}\bm{u},\bm{\zeta})=(\nabla \bm{u},\nabla\bm{\zeta}),\quad  \forall\, \bm{\zeta} \in \bm{H}^1_{0,\mathrm{div}}(\Omega),\nonumber
 \ee
 with domain $D(\bm{S})= \bm{H}^1_{0,\mathrm{div}}(\Omega)\cap\bm{H}^2(\Omega)$ (see e.g., \cite[Chapter III]{S}). 
 The operator $\bm{S}$ is a canonical isomorphism from $\bm{H}^1_{0,\mathrm{div}}(\Omega)$ to $\bm{H}^1_{0,\mathrm{div}}(\Omega)'$. Denote its inverse map by $\bm{S}^{-1}:\bm{H}^1_{0,\mathrm{div}}(\Omega)'\to\bm{H}^1_{0,\mathrm{div}}(\Omega)$. For any $\bm{f}\in \bm{H}^1_{0,\mathrm{div}}(\Omega)'$, there is a unique $\bm{u}=\bm{S}^{-1}\bm{f}\in\bm{H}^1_{0,\mathrm{div}}(\Omega)$ such that
\be
(\nabla\bm{S}^{-1}\bm{f},\nabla \bm{\zeta})=\langle\bm{f},\bm{\zeta}\rangle_{\bm{H}^1_{0,\mathrm{div}}},\quad \forall\, \bm{\zeta} \in \bm{H}^1_{0,\mathrm{div}}(\Omega).\nonumber
\ee
Then we can see that $\|\nabla\bm{S}^{-1}\bm{f}\|=\langle\bm{f},\bm{S}^{-1}\bm{f}\rangle_{\bm{H}^1_{0,\mathrm{div}}}^{\frac{1}{2}}$ is an equivalent norm on $\bm{H}^1_{0,\mathrm{div}}(\Omega)'$ and there exists the chain rule
\be
\langle\bm{f}_{t}(t),\bm{S}^{-1}\bm{f}(t)\rangle_{\bm{H}^1_{0,\mathrm{div}}}=\frac{1}{2}\frac{d}{dt}\|\nabla\bm{S}^{-1}\bm{f}\|^2,\quad   \textrm{a.e.}\ t \in (0,T),\nonumber
\ee
for any $\bm{f}\in H^1(0,T;\bm{H}^1_{0,\mathrm{div}}(\Omega)')$. 
Besides, we recall the following useful result (see e.g.,  \cite[Chapter III, Theorem 2.2.1]{S} and \cite[Appendix B]{GMT}):
\bl \label{stokes}
 \rm Let $d=2,\,3$. For any $\bm{f} \in \bm{L}^2_{\mathrm{div}}(\Omega)$, 
there exists a unique $\bm{u}\in \bm{H}^1_{0,\mathrm{div}}(\Omega)\cap\bm{H}^2(\Omega)$ and $p\in H^1(\Omega)\cap L_0^2(\Omega)$ such that $-\Delta \bm{u}+\nabla p=\bm{f}$ a.e. in $\Omega$, that is, $\bm{u}=\bm{S}^{-1}\bm{f}$. Moreover, 
\begin{align*}
&\|\bm{u}\|_{\bm{H}^2}+\|\nabla p\|\le C\|\bm{f}\|,
\\
& \|p\|\le C \|\bm{f}\|^\frac12\|\nabla \bm{S}^{-1}\bm{f}\|^\frac12,
 \end{align*}
where $C$ is a positive constant that may depend on $d$, $\Omega$ but is independent of $\bm{f}$.
 \el

\subsection{Main results}
\noindent 
We make the following hypotheses.
\begin{enumerate}
\item[(H1)]\label{seta} The viscosity $\eta \in C^{1}(\mathbb{R})$ and satisfies\\
\be 
\eta_{*} \leq \eta(r)\leq \eta^*,\quad |\eta'(r)|\leq \eta_{0},\quad \forall\, r \in \mathbb{R},\nonumber
\ee 
where $\eta_{*}$, $\eta^*$ and $\eta_{0}$ are some positive constants.
\item[(H2)]\label{as} The singular potential $\varPsi$ belongs to the class of functions $C[-1,1]\cap C^{2}(-1,1)$ and can be written into the following form 
\begin{equation}
\varPsi(r)=\varPsi_{0}(r)-\frac{\theta_{0}}{2}r^2,\nonumber
\end{equation}
such that 
\begin{equation}
\lim_{r\to \pm 1} \varPsi_{0}'(r)=\pm \infty ,\quad \text{and}\ \  \varPsi_{0}''(r)\ge \theta,\quad \forall\, r\in (-1,1),\nonumber
\end{equation}
where $\theta$ is a strictly positive constant and $\theta_0\in \mathbb{R}$. In addition, there exists $\epsilon_0\in(0,1)$ such that $\varPsi_{0}''$ is nondecreasing in $[1-\epsilon_0,1)$ and nonincreasing in $(-1,-1+\epsilon_0]$. 
Finally, we make the extension $\varPsi_{0}(r)=+\infty$ for any $r\notin[-1,1]$.
\item[(H3)] The function $h\in C^1(\mathbb{R})\cap L^\infty(\mathbb{R})$ and $S\in L^2(0,T; L^2(\Omega))$.
\item[(H4)]\label{sco} The coefficients $A,\ B,\ \mathcal{C},\ \chi,\ \alpha,\ c_0$ are  prescribed  constants and satisfy
\be 
A>0,\ \ B>0,\ \ \mathcal{C}\in\mathbb{R},\ \ \chi \in \mathbb{R},\ \ \alpha\geq 0,\ \ c_0\in(-1,1). \nonumber
\ee
\end{enumerate}
\begin{remark}
The logarithmic potential \eqref{pot} fulfills the assumption \rm{(H2)}.  
As indicated in \cite[Remark 2.1]{GMT}, one can easily extend the linear viscosity function \eqref{vis} to $\mathbb{R}$ in such a way to comply \rm{(H1)}. Indeed, since the singular potential guarantees that the solution $\varphi\in [-1,1]$, the value of $\eta$ outside of $[-1,1]$ is not important and can be chosen in a good manner as in \rm{(H1)}. Besides, it is possible to consider other physically relevant viscosities like (e.g., \cite{GLLW})
\be
	\eta(r)=\frac{\eta_{1}\eta_{2}}{\eta_{1}(\frac{1-r}{2})+\eta_{2}(\frac{1+r}{2})},\quad  \textrm{or}\ \ \eta(r)=\eta_{1}e^{(\log(\frac{\eta_{1}}{\eta_{2}})(\frac{1-r}{2}))},\quad \forall\, r\in[-1,1],\nonumber
	\ee
	where $\eta_1$ and $\eta_2$ are the viscosities of fluid $1$ and fluid $2$, respectively.
\end{remark}

Next, we introduce the definition of weak solution.
\bd \label{maind}\rm  
Let $d =2,3$ and $T \in (0,+\infty)$. Suppose that the initial data satisfy $\bm {v}_{0} \in \bm {L}^2_{\mathrm{div}}(\Omega)$, $\varphi_{0}\in H^1(\Omega)$, $\sigma_{0}\in L^2(\Omega)$ with $\|  \varphi_{0} \|_{L^{\infty}} \le 1$ and 
$|\overline{\varphi}_{0}|<1$. A quadruple $(\bm{v},\varphi,\mu,\sigma)$ satisfying the following properties
\begin{align}
&\bm{v} \in L^{\infty}(0,T;\bm{L}^2_{\mathrm{div}}(\Omega)) \cap L^{2}(0,T;\bm{H}^1_{0,\mathrm{div}}(\Omega))\cap  W^{1,\frac{4}{d}}(0,T;\bm{H}^1_{0,\mathrm{div}}(\Omega)')\notag,\\
&\varphi \in L^{\infty}(0,T;H^1(\Omega))\cap L^{4}(0,T;H^2_{N}(\Omega))\cap L^2(0,T;W^{2,q}(\Omega)) \cap H^{1}(0,T;H^1(\Omega)'),\notag \\
&\mu \in   L^{2}(0,T;H^1(\Omega)),\notag \\
&\sigma  \in L^{\infty}(0,T;L^2(\Omega))\cap L^{2}(0,T;H^1(\Omega)) \cap W^{1,\frac{4}{d}}(0,T;H^1(\Omega)'),\notag\\
&\varphi\in L^{\infty}(\Omega\times (0,T))\ \textrm{and}\ \ |\varphi(x,t)|<1\ \ \textrm{a.e.\ in}\ \Omega\times(0,T),\notag
\end{align} 
where $q\geq 2$ if $d=2$ and $q\in [2,6]$ if $d=3$, is a weak solution to problem \eqref{f3.c}--\eqref{ini0} on $[0,T]$, if
\begin{subequations}
	\begin{alignat}{3}
	&\left \langle\partial_t  \bm{ v},\bm{\zeta}\right \rangle_{\bm{H}^1_{0,\mathrm{div}}}+(( \bm{ v} \cdot \nabla)  \bm {v},\bm{ \zeta})+(  2\eta(\varphi) D\bm{v},D\bm{ \zeta}) \notag \\
	&\quad=((\mu+\chi \sigma)\nabla \varphi,\bm {\zeta}),\quad &\quad \textrm{a.e.\ in}\ (0.T), &\label{test3.c} \\
	&\left \langle \partial_t \varphi,\xi\right \rangle_{H^1}+((\bm{v} \cdot \nabla) \varphi,\xi)=- (\nabla \mu,\nabla \xi)-\alpha(\varphi-c_0,\xi),\quad  &\quad \textrm{a.e.\ in}\ (0.T),&\label{test1.a} \\
	&\ \, \mu=A\varPsi'(\varphi)-B\Delta \varphi-\chi \sigma,\quad   &\quad \textrm{a.e.\ in}\ (0.T),&\label{test4.d}\\
	&\left \langle\partial_t \sigma,\xi\right \rangle_{H^1}+((\bm{v} \cdot \nabla) \sigma,\xi) + (\nabla \sigma,\nabla \xi)\notag\\
	&\quad = \chi ( \nabla \varphi,\nabla \xi)-(\mathcal{C} h(\varphi) \sigma,\xi) + (S,\xi),\quad   &\quad \textrm{a.e.\ in}\ (0.T),& \label{test2.b} 
	\end{alignat}
\end{subequations}
for all $\bm {\zeta} \in \bm{H}^1_{0,\mathrm{div}}$ and $\xi \in H^1(\Omega)$. Moreover, the initial conditions are fulfilled 
$$\bm {v}(0)=\bm{v}_{0},\quad \varphi(0)=\varphi_{0,},\quad  \sigma(0)=\sigma_{0}.$$
\ed
\begin{remark} 
The initial data are attained in the following sense (see e.g., \cite{B}): from regularity properties of the weak solution and the Sobolev  embedding theorem, we have 
\begin{align*}
&\bm{v} \in C_w([0,T];\bm{L}^2_{\mathrm{div}}(\Omega))\ \ \text{if}\ d=3,\quad \bm{v} \in C([0,T];\bm{L}^2_{\mathrm{div}}(\Omega))\ \ \text{if}\ d=2,\\
& \varphi\in C_w([0,T]; H^1(\Omega)),\\
& \sigma \in C_w([0,T];L^2(\Omega))\ \ \text{if}\ d=3,\quad \sigma \in C([0,T];L^2(\Omega))\ \ \text{if}\ d=2.
\end{align*}
\end{remark}
	
We	 are now in a position to state the main results of this paper.
\bt(Existence of global weak solutions). \label{main}\rm Let $d=2,3$, $T>0$. Suppose that the hypotheses (H1)--(H4) are satisfied, then for any initial data satisfying $\bm {v}_{0} \in \bm {L}^2_{\mathrm{div}}(\Omega)$, $\varphi_{0}\in H^1(\Omega)$, $\sigma_{0}\in L^2(\Omega)$ with $\|  \varphi_{0} \|_{L^{\infty}} \le 1$ and 
$|\overline{\varphi}_{0}|<1$, the initial boundary value problem \eqref{f3.c}--\eqref{ini0} admits at least one global weak solution $(\bm{v},\varphi,\mu,\sigma)$ on $[0,T]$ in the sense of Definition \ref{maind}.
\et

\bt(Continuous dependence estimate with respect to initial data in 2D). \label{the2}
\rm Let $d=2$. Consider two groups of initial data satisfying  $(\bm{v}_{0i},\varphi_{0i},\sigma_{0i})\in\bm{L}^2_{\mathrm{div}}(\Omega)\times H^1(\Omega)\times L^2(\Omega)$ with $\left \|  \varphi_{0i}\right \|_{L^{\infty}} \le 1$, $|\overline{\varphi}_{0i}|<1$, $i=1,\, 2$, and $\overline{\varphi}_{01},\overline{\varphi}_{02}\in (-1,1)$.
The global weak solutions $(\bm{v}_{1},\varphi_{1},\sigma_{1})$,   $(\bm{v}_{2},\varphi_{2 },\sigma_{2})$ to problem \eqref{f3.c}--\eqref{ini0} on $[0,T]$ with initial data $(\bm{v}_{0i},\varphi_{0i},\sigma_{0i})$, $i=1,\, 2$ (and the same source term $S_1=S_2$), 
satisfy the following continuous dependence estimate:
\be
W(t)\le C\left(\frac{W(0)}{C}\right )^{\exp\big(-C\int_{0}^{t}Z(s)\, ds\big)},\quad \forall\, t\in [0,T],
\notag
\ee
where
\begin{align*}
W(t)&= \frac{1}{2}\|\nabla\bm{S}^{-1}[\bm{v}_1(t)-\bm{v}_2(t)]\|^2+\frac{1}{2}\|\varphi_1(t)-\varphi_2(t) \|_{(H^1)'}^2  +\frac{1}{2}\|\sigma_1(t)-\sigma_2(t) \|_{(H^1)'}^2\\
&\quad +|\overline{\varphi}_1(t)-\overline{\varphi}_2(t)|,\\
Z(t)&= \|\nabla \bm{v}_{1}(t)\|^2+\|\nabla \bm{v}_{2}(t)\|^2+\|\varphi_{1}(t)\|_{W^{2,3}}^2+\|\varphi_{2}(t)\|_{W^{2,3}}^2\notag\\
&\quad +\|\varphi_{1}(t)\|_{H^{2}}^4+\|\varPsi'(\varphi_{1})\|_{L^1}+\|\varPsi'(\varphi_{2})\|_{L^1}+\|\sigma_2(t)\|_{H^1}^2+1,
\end{align*}
and $C>0$ is a constant depending on the initial data, $\Omega$ and coefficients of the system.
\et
\begin{corollary}(Uniqueness of weak solutions in 2D).\label{uniq}
\rm Let $d=2$. The global weak solution $(\bm{v},\varphi,\mu,\sigma)$ to problem \eqref{f3.c}--\eqref{ini0} obtained in Theorem \ref{main} is unique.  
\end{corollary}
%

\section{Existence of Global Weak Solutions}\label{ws}
\setcounter{equation}{0}
In this section, we prove Theorem \ref{main} on the existence of global weak solutions to problem \eqref{f3.c}--\eqref{ini0}. The proof relies on a suitable semi-Galerkin scheme. Roughly speaking, the procedure consists of the following steps: first, given a smooth velocity field $\bm{u}^m$, we solve the Cahn--Hilliard equation for $\varphi$ and the reaction-diffusion equation for $\sigma$ with convection terms; second, using the solutions $(\varphi^m, \sigma^m)$ obtained in the previous step, we solve a finite dimensional approximation of the Navier--Stokes equation for $\bm{v}$ with an external force term (given by $\varphi^m$); third, we apply Shauder's fixed point theorem to find a fixed point $(\bm{v}^m, \varphi^m, \sigma^m)$; finally, we derive uniform estimates with respect to $m$ and pass to the limit as $m \to \infty$.

In the subsequent proof, we will use the following modified Gronwall's lemma derived in \cite[Lemma 3.1]{GL17e}.
\bl\label{GronL}
Let $\alpha$, $\beta$, $u$ and $v$ be real-valued functions defined on $[0,T]$. Assume that $\alpha$ is integrable, $\beta$ is non-negative and continuous, $u$ is continuous, $v$ is non-negative and integrable. Suppose $u$ and $v$ satisfy the integral inequality
$$u(s)+\int_0^s v(t)\,dt\leq \alpha(s)+\int_0^s \beta(t)u(t) dt,\quad \forall\,s\in [0,T],
$$
then
\begin{align}
u(s)+\int_0^s v(t)\,dt\leq \alpha(s)+\int_0^s\alpha(t)\beta(t) \exp\left(\int_t^s\beta(r)\,dr\right)dt.\label{gron}
\end{align}
\el

\subsection{Semi-Galerkin Scheme}

 Let the family $\{\bm{y}_{k}(x)\}_{k=1}^{\infty}$ be a basis of the Hilbert space $\bm{H}^1_{0,\mathrm{div}}(\Omega)$, which is given by eigenfunctions of the Stokes problem
\be
(\nabla \bm{y}_{k},\nabla \bm{w})=\lambda_{k}(\bm{y}_{k},\bm{w}),\quad  \forall\, \bm{w} \in {\bm{H}^1_{0,\mathrm{div}}(\Omega)},\quad \textrm{with}\ \|\bm{y}_{k}\|=1,
\ee
where $\lambda_{k}$ is the eigenvalue corresponding to $\bm{y}_{k}$.  It is well-known that $0<\lambda_{1}< \lambda_{2}<... $ is an unbounded monotonically increasing sequence, $\{\bm{y}_{k}(x)\}_{k=1}^{\infty}$ forms a complete orthonormal basis in $\bm{L}^2_{\mathrm{div}}(\Omega)$ and it is also orthogonal in $\bm{H}^1_{0,\mathrm{div}}(\Omega)$. By the elliptic regularity theory, we have $\bm{y}_{k}(x) \in C^{\infty} $ for all $k\in \mathbb{N}$. 
For every $m\in \mathbb{N}$, we denote the finite-dimensional subspace of $\bm{H}^1_{0,\mathrm{div}}(\Omega)$ by 
$$ \bm{H}_{m}:=\textrm{span} \{\bm{y}_{1}(x) ,...,\bm{y}_{m}(x)\}.$$
Moreover, we use $\bm{P}_{\bm{H}_m}$ for the corresponding orthogonal projections from $\bm{L}^2_{\mathrm{div}}(\Omega)$ onto $\bm{H}_m$.

For every $m\in \mathbb{N}$ and arbitrary $T>0$, we consider the following approximate problem: looking for functions
\be
\bm{v}^{m}(x,t):=\sum_{i=1}^{m}a_{i}^{m}(t)\bm{y}_{i}(x),\label{appovm}
\ee
and $(\varphi^{m},\, \mu^{m},\, \sigma^{m})$ satisfying
\begin{numcases}
{\textrm{(P1)}}
(\partial_t  \bm{ v}^{m},\bm{w})+(( \bm{ v}^{m} \cdot \nabla)  \bm {v}^{m},\bm{ w})+(  2\eta(\varphi^{m}) D\bm{v}^{m},D\bm{w}) \notag\\
\quad=((\mu^{m}+\chi \sigma^{m})\nabla \varphi^{m},\bm {w}),\qquad \qquad \qquad \qquad \quad \,
\text{a.e. in}\ (0,T), \label{ptest3.c}\\
\left \langle\partial_t \varphi^m,\xi\right \rangle_{{H}^1}+((\bm{v}^{m} \cdot \nabla) \varphi^{m},\xi)\notag\\
\quad =-(\nabla \mu^{m},\nabla\xi)-\alpha(\varphi^m-c_0,\xi),\quad\qquad \qquad \quad  \text{a.e. in}\ (0,T),\label{pag4.d} \\
\mu^{m}=A\varPsi'(\varphi^{m})-B\Delta \varphi^{m}-\chi \sigma^{m}, \qquad \qquad\quad \quad \ \ \rm{a.e.\ in}\ \Omega\times(0,T),\label{pag1.a}\\
\left \langle\partial_t \sigma^m,\xi\right \rangle_{{H}^1}+((\bm{v}^{m} \cdot \nabla) \sigma^{m},\xi)+(\nabla \sigma^m,\nabla \xi)\notag\\
\quad = \chi ( \nabla \varphi^m,\nabla \xi)-(\mathcal{C} h(\varphi^m) \sigma^m,\xi) + (S,\xi),\qquad\,   \text{a.e. in}\ (0,T),\label{pag2.b}\\
\bm{v}^{m}(0)=\bm{P}_{\bm{H}_{m}} \bm{v}_{0},\quad \varphi^{m}(0)=\varphi_{0,},\quad  \sigma^{m}(0)=\sigma_{0},\quad \textrm{in}\ \Omega,
\end{numcases}
for all $\bm{w} \in \bm{H}_{m}$, $\xi \in H^1(\Omega)$.

The following proposition yields that the approximate problem (P1) admits a unique weak solution.  
\bp\label{p1}\rm Let $d=2,3$. We assume that the hypotheses (H1)--(H4) are satisfied, and the initial data satisfy $\bm {v}_{0} \in \bm {L}^2_{\mathrm{div}}(\Omega)$, $\varphi_{0}\in H^1(\Omega)$, $\sigma_{0}\in L^2(\Omega)$ with $\|  \varphi_{0} \|_{L^{\infty}} \le 1$ and $|\overline{\varphi}_{0}|<1$. 
For every integer $m>0$, there exists a time $T_{m}>0$ depending on $\bm{v}_{0},\ \varphi_{0},\ \sigma_{0}$, $\Omega$, $m$ and coefficients of the system such that problem (P1) admits a unique weak solution $(\bm{v}^{m},\varphi^{m},\mu^{m},\sigma^{m})$ on $[0,T_{m}]$ satisfying
\begin{align}
&\bm{v}^{m} \in H^1(0,T_{m};\bm{H}_m(\Omega))\notag,\\
&\varphi^{m} \in L^{\infty}(0,T_{m};H^1(\Omega))\cap L^{2}(0,T_{m};H^2_{N}(\Omega)) \cap H^{1}(0,T_{m};H^1(\Omega)'),\notag \\
&\mu^{m} \in   L^{2}(0,T_{m};H^1(\Omega)),\notag \\
&\sigma^{m}   \in L^{\infty}(0,T_{m};L^2(\Omega))\cap L^{2}(0,T_{m};H^1(\Omega)) \cap H^1(0,T_{m};H^1(\Omega)').\notag \\
&\varphi^m \in L^{\infty}(\Omega\times(0,T))\ \textrm{and}\ |\varphi|<1\ \textrm{a.e.\ in}\ \Omega\times(0,T).\notag
\end{align}
\ep
\begin{remark}
In view of \eqref{appovm}, the approximate solution $\bm{v}^m$ is indeed smooth in space.
\end{remark}
\noindent\textbf{Proof.} 
The proof of Proposition \ref{p1} consists of several steps. \medskip

\textbf{Step 1}. Let $T>0$ and  $\widetilde{M}\geq 2\|\bm{v}_0\|^2+1$ to be determined later. Consider an arbitrary given function
\be
\bm{u}^{m}=\sum_{i=1}^{m}a_{i}^{m}(t)\bm{y}_{i}(x)\in C([0,T];\bm{H}_{m}),\nonumber
\ee
that satisfies
\be
a_{i}^{m}(0)=(\bm{v}_{0},\bm{y}_{i}),\quad \sup_{t\in [0,T]}  \sum_{i=1}^{m}|a_{i}^{m}(t)|^{2}\le \widetilde{M}.\nonumber
\ee
Namely, $\bm{u}^{m}(0)= \bm{P}_{\bm{H}_{m}} \bm{v}_{0}$ and 
$\sup_{t\in [0,T]}\|\bm{u}^{m}(t)\|^2\le \widetilde{M}$. Besides, since $\bm{u}^m$ is indeed finite dimensional, we have the inverse inequality
$$
  \sup_{t\in [0,T]}\|\bm{u}^{m}(t)\|^2_{\bm{L}^{\infty}}\le m \widetilde{M} \max_{1\le i\le m}\|\bm{y}_{i}\|^2_{\bm{L}^{\infty}}\le C_{m}\widetilde{M}, 
$$
where $C_m>0$ is a constant that depends on $m$.

With the given velocity vector $\bm{u}^{m}$, we first consider the following auxiliary system:
\begin{subequations}
	\begin{alignat}{3}
	&\left \langle\partial_t  \varphi^{m},\xi\right \rangle_{{H}^1}+((\bm{u}^{m} \cdot \nabla) \varphi^{m},\xi)\notag\\
	&\quad =-(\nabla \mu^{m},\nabla\xi)-\alpha(\varphi^m-c_0,\xi),\qquad \qquad \qquad \qquad \quad  \textrm{a.e.\ in}\ (0,T),\label{g1.a} \\
	&\mu^{m}=A\varPsi'(\varphi^{m})-B\Delta \varphi^{m}-\chi \sigma^{m}, \qquad \qquad \qquad \qquad \qquad\rm{a.e.\ in}\ \Omega\times(0,T),\label{g4.d} \\
	&\left \langle\partial_t  \sigma^{m},\xi\right \rangle_{{H}^1}+((\bm{u}^{m} \cdot \nabla) \sigma^{m},\xi)+(\nabla \sigma^m,\nabla \xi)\notag\\
&\quad = \chi ( \nabla \varphi^m,\nabla \xi)-(\mathcal{C} h(\varphi^m) \sigma^m,\xi) + (S,\xi),\quad\, \qquad \qquad \text{a.e. in}\ (0,T), \label{g2.b}\\
	&\varphi^{m}(0)=\varphi_{0},\quad  \sigma^{m}(0)=\sigma_{0},\ \ \qquad \qquad \qquad \qquad \qquad \qquad\textrm{in}\ \Omega,  \label{g6.h}
	\end{alignat}
\end{subequations}
for all $\xi \in H^1(\Omega)$. Then we have 
\bl\label{fp}
\rm  Assume that $\bm{u}^{m}\in C([0,T];\bm{H}_{m})$ as described as above. For any initial data $\varphi_{0}\in H^1(\Omega)$, $\sigma_{0}\in L^2(\Omega)$ with $\|  \varphi_{0} \|_{L^{\infty}} \le 1$ and $|\overline{\varphi}_{0}|<1$ given in the statement of Proposition \ref{p1}, the auxiliary problem \eqref{g1.a}--\eqref{g6.h}  admits a unique weak solution $(\varphi^{m},\mu^{m},\sigma^{m})$  satisfying
\begin{align}
&\varphi^{m} \in C_w([0,T];H^1(\Omega))\cap L^{2}(0,T;H^2_{N}(\Omega)) \cap H^{1}(0,T;H^1(\Omega)'),\notag \\
&\mu^{m} \in   L^{2}(0,T;H^1(\Omega)),\notag \\
&\sigma^{m}   \in C([0,T];L^2(\Omega))\cap L^{2}(0,T;H^1(\Omega)) \cap H^1(0,T;H^1(\Omega)'),\notag \\
&\varphi^m\in   L^{\infty}(\Omega\times(0,T))\ \ \textrm{and}\ \ |\varphi^{m}|<1\quad \textrm{a.e.\ in}\ \Omega\times(0,T).\notag
\end{align} 
\el
The proof of Lemma \ref{fp} consists of several steps and is postponed to the Appendix. Next, we introduce the following Banach space 
\begin{align}
X&=\big(L^{\infty}(0,T;H^1(\Omega)) \cap  L^{2}(0,T;H^2_{N}(\Omega))\big)\notag \\
&\quad\times \big(L^{\infty}(0,T;L^2(\Omega))\cap  L^{2}(0,T;H^1(\Omega))\big),\notag
\end{align}
and the mapping 
\begin{align}
\Phi^m_{1}:\  C([0,T];\bm{H}_{m})\ &\to\ \  X, \notag\\
\bm{u}^{m}\ &\to\ \ (\varphi^{m},\sigma^{m}).\notag 
\end{align}
Thanks to Lemma \ref{fp}, the mapping $\Phi^m_{1}$ is well-defined, and $\Phi^m_{1}(\bm{u}^{m})=(\varphi^{m},\sigma^{m})$ is bounded from $C([0,T];\bm{H}_{m})$ to $X$. 

We proceed to show that $\Phi^m_{1}$ is also continuous. To this end, let $\bm{u}^{m}_{1}$ and $\bm{u}^{m}_{2} \in C([0,T];\bm{H}_{m})$ be two given vectors with the same initial value as above. Then $(\varphi^{m}_{i},\sigma^{m}_{i}) =\Phi^m_{1}(\bm{u}^{m}_{i})$, $i=1,\, 2$ are the two weak solutions to problem \eqref{g1.a}--\eqref{g6.h} given by Lemma \ref{fp} (subject to the same initial data $(\varphi_0,\sigma_0)$ and source term $S$), with the corresponding chemical potentials $\mu^m_i$ being given by \eqref{g4.d}. 
We denote the differences by 
$$
\bm{u}^m= \bm{u}^m_1-\bm{u}^m_2,\quad (\varphi^{m},\mu^{m},\sigma^{m})=(\varphi_{1}^{m}-\varphi_{2}^{m},\mu_{1}^{m}-\mu_{2}^{m},\sigma_{1}^{m}-\sigma_{2}^{m}).
$$
Then it holds
\begin{subequations}
	\begin{alignat}{3}
	&\left \langle\partial_t  \varphi^{m},\xi\right \rangle_{ {H}^1}+((\bm{u}^{m}_{1} \cdot \nabla) \varphi^{m},\xi)+((\bm{u}^{m}\cdot \nabla) \varphi^{m}_{2},\xi)\notag\\
	&\quad =- (\nabla \mu^{m},\nabla \xi)-\alpha(\varphi^m,\xi),\label{atest11.a} \\
	&(\mu^{m},\xi)=A(\varPsi'(\varphi^{m}_{1})-\varPsi'(\varphi^{m}_{2}),\xi)+B(\nabla \varphi^{m},\nabla \xi)-(\chi \sigma^{m},\xi) ,\label{atest44.d}\\
	&\left \langle\partial_t  \sigma^{m},\xi\right \rangle_{{H}^1}+((\bm{u}^{m}_{1} \cdot \nabla) \sigma^{m},\xi)+((\bm{u}^{m} \cdot \nabla) \sigma^{m}_{2},\xi)+(\nabla \sigma^m,\nabla\xi)
	\notag\\
	&\quad =\chi(  \nabla \varphi^{m},\nabla \xi)
	-\mathcal{C}( h(\varphi^m_1) \sigma^m_1- h(\varphi^m_2)\sigma^m_2,\xi), \label{atest22.b}\\
	& \varphi^{m}(0)=0,\quad  \sigma^{m}(0)=0.
	\end{alignat}
\end{subequations}
After integration by parts, and using the fact $\mathrm{div}\,\bm{u}^m=0$, \eqref{atest11.a} can be rewritten as
\begin{equation}
\left \langle\partial_t  \varphi^{m},\xi\right \rangle_{{(H^1)}', {H}^1}-(\varphi^{m}\bm{u}^{m}_{1},\nabla\xi)-(\varphi^{m}_{2}\bm{u}^{m},\nabla\xi)=- (\nabla \mu^{m},\nabla \xi)-\alpha(\varphi^m,\xi).\label{atest111.a}
\end{equation}
From Lemma \ref{fp}, we have the estimates
\begin{equation}
 \|\varphi^{m}_{i}(t)\|_{H^1} \le C,\quad  \|\varphi^{m}_{i}(t)\|_{L^{\infty}} \le 1,\quad \forall\, t \in[0,T],\quad i=1,2.\nonumber
\end{equation}
Besides, taking $\xi=1$ in \eqref{g1.a}, we easily get
$$\overline{\varphi}_i^m(t)=c_0+e^{-\alpha t}(\overline{\varphi}_0-c_0),\quad i=1,2,
$$
which implies
\be
\overline{\varphi}^{m}(t)=\overline{\varphi}^{m}_{1}(t)-\overline{\varphi}^{m}_{2}(t)=0,\quad \forall\, t \in[0,T].\notag
\ee
Now choosing the test function $\xi=\mathcal{N}\varphi^{m}$ in \eqref{atest111.a}, we obtain
\begin{equation}
\frac{1}{2}\frac{d}{dt}\|\varphi^{m}\|_{V_0'}^2+\alpha\|\varphi^m\|_{V_0'}+(\mu^{m},\varphi^{m})=I_{1}+I_{2},\label{vhmna}
\end{equation}
where 
$$I_{1}=(\varphi^{m}\bm{u}^{m}_{1},\nabla \mathcal{N}\varphi^{m}),\ \ I_{2}=(\varphi^{m}_{2}\bm{u}^{m},\nabla \mathcal{N}\varphi^{m}).$$
From the assumption (H2) on $\varPsi$, we deduce that
\begin{align}
(\mu^{m},\varphi^{m})&=A(\varPsi'(\varphi^{m}_{1})-\varPsi'(\varphi^{m}_{2}),\varphi^{m})+B(\nabla \varphi^{m},\nabla \varphi^{m})-(\chi \sigma^{m},\varphi^{m}) \notag\\
&\ge B\|\nabla \varphi^{m}\|^2-C_1\|\varphi^{m}\|^2-\frac{1}{4}\|\sigma^{m}\|^2,\notag
\end{align}
where $C_1$ is a constant depending on $\theta$, $\theta_0$ and $\chi$. By the definition of $\mathcal{N}$ and Young's inequality, we have
\begin{align}
C_1\|\varphi^{m}\|^2&=C_1(\nabla A_{0}^{-1}\varphi^{m},\nabla \varphi^{m})\le\frac{B}{8}\|\nabla \varphi^{m}\|^2+C_2\|\varphi^{m}\|_{V_0'}^2.\notag
\end{align}
Besides, we have the following estimates
\begin{align}
I_{1}
&\le \|\varphi^m\|_{L^6}\|\bm{u}^m_1\|_{\bm{L}^3}\|\varphi^m\|_{V_0'}\nonumber\\
&\le\frac{B}{8}\|\nabla\varphi^{m}\|^2+C\|\bm{u}^m_{1}\|_{\bm{L}^3}^2\|\varphi^{m}\|_{V_0'}^2\notag\\
&\le\frac{B}{8}\|\nabla\varphi^{m}\|^2+C\|\varphi^{m}\|_{V_0'}^2,\notag
\end{align}
and
\begin{align}
I_{2}
&\le \|\varphi^m_2\|_{L^6}\|\bm{u}^m\|_{\bm{L}^3}\|\varphi^m\|_{V_0'}
\le\frac{1}{4}\|\bm{u}^{m}\|^2+C\|\varphi^{m}\|_{V_0'}^2,\notag
\end{align}
where we have used the inequality $\|\bm{u}^m\|_{\bm{L}^3}\leq |\Omega|^\frac13\|\bm{u}^m\|_{\bm{L}^\infty}\leq  C_m\|\bm{u}^m\|$, since it is indeed finite dimensional.
Collecting the above estimates, we obtain that 
\begin{equation}
\frac{1}{2}\frac{d}{dt}\|\varphi^{m}\|_{V_0'}^2+\frac{3}{4}B\|\nabla \varphi^{m}\|^2\le C\|\varphi^{m}\|_{V_0'}^2+\frac{1}{4}\|\sigma^{m}\|^2+\frac{1}{4}\|\bm{u}^{m}\|^2.
\label{astarphi}
\end{equation}
%
Next, taking $\xi=\mathcal{N}_1\sigma^{m}$ in \eqref{atest22.b}, we get
\begin{align}
&\frac{1}{2}\frac{d}{dt} \|\sigma^m\|_{(H^1)'}^2+\|\sigma^m\|^2\nonumber\\
&\quad = \|\sigma^m\|_{(H^1)'}^2
+(\bm{u}^{m}_{1} \sigma^{m}, \nabla \mathcal{N}_1\sigma^{m})
+(\bm{u}^{m} \sigma^{m}_2, \nabla \mathcal{N}_1\sigma^{m})\nonumber\\
&\qquad +\chi(\nabla \varphi^m,\nabla \mathcal{N}_1\sigma^{m})
-\mathcal{C}(h(\varphi^m_1)\sigma^m,\mathcal{N}_1\sigma^{m})-\mathcal{C}((h(\varphi^m_1)-h(\varphi^m_2))\sigma^m_2,\mathcal{N}_1\sigma^{m})\nonumber\\
&\quad := \|\sigma^m\|_{(H^1)'}^2+\sum_{j=3}^7 I_j,\label{diffsig}
\end{align}
where
\begin{align*}
I_3&=(\bm{u}^{m}_{1} \sigma^{m}, \nabla \mathcal{N}_1\sigma^{m}),\qquad\quad  
I_4=(\bm{u}^{m} \sigma^{m}_2, \nabla \mathcal{N}_1\sigma^{m}),\\
I_{5}&=\chi(\nabla \varphi^m,\nabla \mathcal{N}_1\sigma^{m}),\quad \qquad
I_{6}=-\mathcal{C}(h(\varphi^m_1)\sigma^m,\mathcal{N}_1\sigma^{m}),\\
I_{7}&=-\mathcal{C}((h(\varphi^m_1)-h(\varphi^m_2))\sigma^m_2,\mathcal{N}_1\sigma^{m}).\quad \
\end{align*}
The right-hand side of \eqref{diffsig} can be estimated as follows: 
\begin{align}
I_3
&\le C\|\bm{u}^{m}_{1}\|_{\bm{L}^\infty}\|\sigma^{m}\|\|\nabla\mathcal{N}_1\sigma^{m}\|\notag\\
&\le C\|\bm{u}^{m}_{1}\|\|\sigma^{m}\|
\|\sigma^{m}\|_{(H^1)'}\nonumber\\
&\leq \frac{1}{8}\|\sigma^m\|^2+C \|\sigma^{m}\|_{(H^1)'}^2,\notag
\end{align}
and similarly,
\begin{align}
I_4
&\le C\|\sigma^{m}_{2}\|\|\bm{u}^{m}\|_{\bm{L}^\infty}\|\nabla \mathcal{N}_1\sigma^{m}\|\notag\\
&\le C\|\bm{u}^{m}\|\|\sigma^{m}\|_{(H^1)'}\notag\\
& \le \frac{1}{4} \|\bm{u}^{m}\|^2+C \|\sigma^{m}\|_{(H^1)'}^2,
\notag
\end{align}
where we again used the inverse inequality $\|\bm{u}^m\|_{\bm{L}^\infty}\leq C_m\|\bm{u}^m\|$. For $I_5$, we infer from Young's inequality that
\begin{align}
I_5
&\le |\chi|\|\nabla \varphi^m\|\|\nabla\mathcal{N}_1\sigma^{m}\| \leq \frac{B}{8}\|\nabla \varphi^m\|^2 + C \|\sigma^{m}\|^2_{(H^1)'}.\nonumber
\end{align}
Concerning $I_6$, $I_7$, we have
\begin{align}
I_6
&\le C \|h(\varphi^m_1)\|_{L^\infty}\|\sigma^m\|\|\mathcal{N}_1\sigma^{m}\|
\leq \frac{1}{8}\|\sigma^m\|^2+C \|\sigma^{m}\|_{(H^1)'}^2,\nonumber
\end{align}
and
\begin{align}
I_7
&\le C\|\mathcal{N}_{1}^\frac12((h(\varphi^m_1)-h(\varphi^m_2))\sigma^m_2)\| \|\mathcal{N}_{1}^\frac12 \sigma^m\|\nonumber \\
&\le C\|(h(\varphi^m_1)-h(\varphi^m_2))\sigma^m_2)\|_{L^\frac{6}{5}}\|\sigma^{m}\|_{(H^1)'}\nonumber\\
&\leq C\|h(\varphi^m_1)-h(\varphi^m_2)\|_{L^3}\|\sigma^m_2\|\|\sigma^{m}\|_{(H^1)'}\nonumber\\
&\le \|h'\|_{L^6}\|\varphi^m\|_{L^6}\|\sigma^{m}\|_{(H^1)'}\nonumber\\
&\le  \frac{B}{8}\|\nabla \varphi^m\|^2 + C \|\sigma^{m}\|^2_{(H^1)'}.\nonumber
\end{align}
Hence, we infer from the above estimates and \eqref{diffsig} 
that 
\begin{align}
&\frac{1}{2}\frac{d}{dt} \|\sigma^m\|_{(H^1)'}^2
+\frac{3}{4}\|\sigma^m\|^2
\le \frac{B}{4}\|\nabla \varphi^m\|^2 +C  \|\sigma^m\|_{(H^1)'}^2 + \frac{1}{4} \|\bm{u}^{m}\|^2.
\label{diffsig1}
\end{align}
In view of \eqref{astarphi} and \eqref{diffsig1}, we deduce the following inequality
\begin{align}
&\frac{d}{dt}\Big(\|\varphi^{m}\|_{V_0'}^2+  \|\sigma^m\|_{(H^1)'}^2\Big)+ B\|\nabla \varphi^{m}\|^2+\|\sigma^m\|^2
\nonumber\\
&\quad \le C\Big(\|\varphi^{m}\|_{V_0'}^2+ \|\sigma^m\|_{(H^1)'}^2\Big)+\|\bm{u}^{m}\|^2,
\label{diffps}
\end{align}
where $C>0$ is a constant depending on the initial data, $\Omega$, $\widetilde{M}$, $m$ and coefficients of the system. 
Integrating \eqref{diffps} with respect to time, we infer from Gronwall's lemma (see Lemma \ref{GronL}) that 
\begin{align}
&\|\varphi^{m}(t)\|_{V_0'}^2+  \|\sigma^m(t)\|_{(H^1)'}^2 +\int_0^t \Big(\|\nabla \varphi^{m}(s)\|^2+\|\sigma^m(s)\|^2\Big) ds\notag\\
&\le Cte^{Ct}\sup_{s\in [0,t]}\|\bm{u}^{m}(s)\|^2,\quad \forall\,t\in[0,T].
\label{auphih1}
\end{align}
Therefore, the solution operator $\Phi^m_{1}: C([0,T];\bm{H}_{m})\to X$ is continuous with respect to $\bm{u}^{m}$ in a weaker topology for the space $X$.
\medskip

\textbf{Step 2}. Once the triple $(\varphi^{m},\mu^{m},\sigma^{m})$ is determined as in Step 1,  we turn to look for functions:
\be
\bm{v}^m=\sum_{i=1}^{m}\tilde{a}_{i}^{m}(t)\bm{y}_{i}(x)\nonumber
\ee
that satisfy the following system
\begin{align}
	&\left\langle\partial_t  \bm{v}^m,\bm{w}\right\rangle_{\bm{H}^1_{0,\mathrm{div}}}
	+(( \bm{v}^m \cdot \nabla)  \bm{v}^m,\bm{ w})+(  2\eta(\varphi^{m}) D\bm{v}^m, D\bm{w}),\notag\\
	&\quad=((\mu^{m}+\chi \sigma^{m})\nabla \varphi^{m},\bm {w}),
	\quad\forall\,\bm{ w}\in\bm{H}_{m},
	\label{aatest3.c} 
	\end{align}
subject to the initial condition
\begin{align}
\bm{v}^m(0)=\bm{P}_{\bm{H}_{m}} \bm{v}_{0}.\label{aatest3.cini} 
\end{align}
Problem \eqref{aatest3.c}--\eqref{aatest3.cini} is equivalent to a system consisting of $m$ nonlinear ordinary differential equations for the coefficients $\{\tilde{a}_{i}^{m}\}_{i=1}^{m}$ (by taking $\bm{w}=\bm{y}_i$, $i=1,...,m$). Since $\bm{ w}\in\bm{H}_{m}$, the external force term reads $((\mu^{m}+\chi \sigma^{m})\nabla \varphi^{m},\bm {w})=-(\varphi^m \nabla(\mu^m+\chi\sigma^m),\bm{w})$. It follows from Lemma \ref{fp} that $\varphi^m \nabla(\mu^m+\chi\sigma^m)\in L^2(0,T;\bm{L}^2(\Omega))$, thus the right-hand side of \eqref{aatest3.c} belongs to $L^2(0,T)$. Due to the smoothness assumption of $\eta$ in (H1), then by the classical theory of ODEs,  it is standard to show the existence and uniqueness of local solutions $\tilde{a}_{i}^{m}\in H^1(0,T_m)$, $i=1,...,m$ to the above ODE system on a certain time interval $[0,T_m]\subset[0,T]$. 

Hence, we obtain a unique local solution $
\bm{v}^m\in H^1(0,T_{m};\bm{H}_{m})$ to problem \eqref{aatest3.c}--\eqref{aatest3.cini}. 
Besides, testing \eqref{aatest3.c} with $\bm{v}^{m}$, we have 
\begin{align}
&\frac12 \frac{d}{dt}\|\bm{v}^m\|^2 +\int_{\Omega}2\eta(\varphi^{m} )|D\bm{v}^{m}|^2 \ dx \notag\\
&\quad = \int_{\Omega}(\mu^{m}+\chi \sigma^{m})\nabla  \varphi^{m} \cdot\bm{v}^m dx\notag \\
&\quad  = -\int_{\Omega} \varphi^m \nabla (\mu^{m}+\chi \sigma^{m}) \cdot\bm{v}^m dx\notag \\
&\quad \leq \|\varphi^m\|_{L^\infty}(\|\nabla \mu^m\|+\|\chi \nabla \sigma^m\|)\|\bm{v}^m\|\notag\\
&\quad \le \frac12 \|\bm{v}^m\|^2+  \|\nabla\mu^{m}\|^2+\chi^2\|\nabla\sigma^{m}\|^2.
\label{aenergy}
\end{align}
Integrating the above inequality with respect to time, from Lemma \ref{GronL}, we deduce that 
\begin{align}
&\|\bm{v}^{m}(t)\|^2 + 4\eta_*\int_{0}^{t}\|D\bm{v}^{m}(s)\|^2 \, ds \notag \\
&\quad\le
C e^{t}\left[\|\bm{v}_{0}\|^2+\int_0^t\big(\|\nabla\mu^{m}(s)\|^2+\|\nabla\sigma^{m}(s)\|^2 \big)ds\right],\quad \forall\,t\in[0,T_m],
\label{auvm1}
\end{align}
where the constant $C>0$ may depend on $\chi$ and $\eta_*$. The above estimate enables us to extend the local solution $\bm{v}^m$ to the whole interval $[0,T]$ such that
\be
\bm{v}^m\in H^1(0,T;\bm{H}_{m}).\label{auvm2}
\ee
Going back to the equation \eqref{aatest3.c} and using the fact $\varphi^m \nabla (\mu^{m}+\chi \sigma^{m}) \in L^2(0,T; \bm{L}^2(\Omega))$, we also infer that 
$$\|\partial_{t}\bm{v}^{m}\|_{L^2(0,T;\bm{H}_{m})}\le C(T).$$
In summary, we have 
\bl\label{NSSa}
Given a triple $(\varphi^{m},\mu^{m},\sigma^{m})$ determined by Lemma \ref{fp}, for any initial datum $\bm{v}_0\in \bm{L}^2_{\mathrm{div}}(\Omega)$ as given in Proposition \ref{p1}, problem \eqref{aatest3.c}--\eqref{aatest3.cini} admits a unique solution $
\bm{v}^m\in H^1(0,T;\bm{H}_{m})$.
\el

Thanks to Lemma \ref{NSSa}, we can define the following mapping determined by the solution to problem  \eqref{aatest3.c}--\eqref{aatest3.cini}:
\begin{align*}
\Phi^m_{2} :\quad\quad X\quad&\to\ H^1(0,T_{m};\bm{H}_{m}),\\
(\varphi^{m},\sigma^{m})\ &\to\ \bm{v}^{m}.
\end{align*}
It is obvious from \eqref{auvm1}, \eqref{auvm2} that $\Phi^m_2$ is bounded from $X$ to $H^1(0,T_{m};\bm{H}_{m})$. Below we verify its continuity. To this end, let $(\varphi^{m}_{i},\sigma^{m}_{i})=\Phi^m_1(\bm{u}^m_i)$, $\bm{u}^m_i\in C([0,T];\bm{H}_{m})$, $i=1,2$, with corresponding $\mu^m_i\in L^2(0,T,H^1(\Omega))$ satisfying $\mu^{m}_i=A\varPsi'(\varphi^{m}_i)-B\Delta \varphi^{m}_i-\chi \sigma^{m}_i$, $i=1,2$ (i.e., given by Lemma \ref{fp}). Then we set 
$$\bm{v}^{m}_{i}=\Phi^m_2(\varphi^{m}_{i},\sigma^{m}_{i}),\quad i=1,\ 2,$$ 
and denote the differences by 
 $$\bm{v}^{m}=\bm{v}^{m}_{1}-\bm{v}^{m}_{2},\quad \quad (\varphi^{m},\sigma^{m})=(\varphi_{1}^{m}-\varphi_{2}^{m},\sigma_{1}^{m}-\sigma_{2}^{m}).$$
Taking the difference of \eqref{aatest3.c} for $\bm{v}^{m}_{i}$ and testing the resultant by $\bm{w}=\bm{v}^m$, we get
\begin{align}
&(\partial_t  \bm{v}^{m},\bm{v}^m)
	+(  2\eta(\varphi^{m}_{1}) D\bm{v}^{m} ,D\bm{v}^m)\notag\\
&\quad  
= 	-(( \bm{ v}^{m} \cdot \nabla)  \bm {v}^{m}_{1},\bm{v}^m)-(( \bm{ v}^{m}_{2} \cdot \nabla)  \bm {v}^{m},\bm{v}^m)
-(2(\eta(\varphi^{m}_{1})-\eta(\varphi^{m}_{2})) D\bm{v}^{m}_{2} ,D\bm{v}^m)\notag\\
&\qquad +\big((\mu^{m}_{1}+\chi\sigma^m_1)\nabla\varphi^{m}_{1}-(\mu^{m}_{2}+\chi\sigma^m_2)\nabla\varphi^{m}_{2},\bm{v}^m\big)\notag\\
&\quad	:=\sum_{i=1}^{4}J_{i},
\label{diffvm}
\end{align} 
where 
\begin{align*}
J_1&=-(( \bm{ v}^{m} \cdot \nabla)  \bm {v}^{m}_{1},\bm{v}^m),\\
 J_{2}&=-(( \bm{ v}^{m}_{2} \cdot \nabla)  \bm {v}^{m},\bm{v}^m),\\
J_{3}&=-(2(\eta(\varphi^{m}_{1})-\eta(\varphi^{m}_{2})) D\bm{v}^{m}_{2} ,D\bm{v}^m),\\
J_{4}&=\big((\mu^{m}_{1}+\chi\sigma^m_1)\nabla\varphi^{m}_{1}-(\mu^{m}_{2}+\chi\sigma^m_2)\nabla\varphi^{m}_{2},\bm{v}^m\big).
\end{align*}
The right-hand side of \eqref{diffvm} can be estimated as follows:
\begin{align}
J_{1}
&\le C \|\nabla \bm{v}^{m}_{1}\|_{\bm{L}^\infty} 
\|\bm{v}^m\|^2\leq C \|\bm{v}^m\|^2,\notag\\
J_{2}
&=0,\notag\\
J_{3}
&\le \sup_{s\in[-1,1]}|2\eta'(s)|\|\varphi^m\|\|D \bm{v}^m_2\|_{\bm{L}^\infty}\|D \bm{v}^m\|\notag\\
& \le C\|\varphi^m\|\|D \bm{v}^m\|\notag\\
& \le \|\bm{v}^m\|^2 + C\|\varphi^{m}\|^2\notag\\
& \le \|\bm{v}^m\|^2 + C\|\nabla \varphi^{m}\|^2.\notag
\end{align}
Using the identity 
\begin{align}
\big((\mu^m_i+\chi\sigma^m_i)\nabla\varphi^m_i, \bm{v}^m\big) 
&= \big((A\varPsi'(\varphi^m_i)-B\Delta \varphi^m_i)\nabla \varphi^m_i, \bm{v}^m\big)\nonumber\\
&= \Big(A \nabla \varPsi(\varphi^m_i) - B\,\mathrm{div}(\nabla \varphi^m_i\otimes\nabla \varphi^m_i)+\frac{B}{2}\nabla |\nabla \varphi^m_i|^2,\, \bm{v}^m\Big) \notag\\
&= -B\big(\mathrm{div}(\nabla \varphi^m_i\otimes\nabla \varphi^m_i),\, \bm{v}^m\big) \nonumber\\
&= B\int_\Omega (\nabla \varphi^m_i\otimes\nabla \varphi^m_i) :\nabla \bm{v}^m dx,\quad i=1,2,\label{eqva1}
\end{align}
we can estimate $J_4$ as follows
\begin{align*}
J_{4}
&= B \int_\Omega [(\nabla \varphi^m_1+\nabla \varphi^m_2)\otimes\nabla \varphi^m] :\nabla \bm{v}^m dx\notag\\
&\le B (\|\nabla \varphi^m_1\|+\|\nabla \varphi^m_2\|)\|\nabla \varphi^m\|\|\nabla \bm{v}^m\|_{\bm{L}^\infty}\notag\\
&\le C\|\nabla \varphi^m\|\|\bm{v}^m\|\notag\\
&\le \|\bm{v}^m\|^2+C \|\nabla \varphi^m\|^2.
\end{align*}
In the above estimates, we always use the fact that $\bm{v}^m_i$ are indeed finite dimensional and the higher-order norms in space can be controlled by their $\bm{L}^2$-norm with a constant that depends on $m$. 

From the boundedness of $\eta$ and above estimates, we deduce from \eqref{diffvm} that 
\begin{align}
\frac12\frac{d}{dt}\|\bm{v}^m\|^2+ 2\eta_*\|D\bm{v}^m\|^2\leq C\|\bm{v}^m\|^2+ C\|\nabla \varphi^m\|^2.
\label{diffvm1}
\end{align}
Again from Lemma \ref{GronL}, we get (keeping in mind that $\bm{v}^m(0)=\bm{0}$)
\be
\|\bm{v}^m(t)\|^2+4\eta_*\int_0^t \|D\bm{v}^m(s)\|^2ds
\le Ce^{Ct} \int_0^t\|\nabla \varphi^m(s)\|^2\,ds,\quad \forall\, t\in [0,T].
\label{aum}
\ee
Therefore, the mapping $\Phi_{2}: X\to H^1(0,T;\bm{H}_{m})$ is continuous with respect to $\varphi^{m}$ (and independent of $\sigma^m$ indeed) in a weaker topology (i.e., w.r.t. $C([0,T];\bm{H}_{m})$).
\medskip

\textbf{Step 3.} We now define the mapping 
\begin{align*}
\Phi^m:=\Phi_{2}^m\circ\Phi_{1}^m :C([0,T];\bm{H}_{m}) \ &\ \to H^1(0,T;\bm{H}_{m}),\\
\bm{u}^m\ &\ \to\bm{v}^{m}.
\end{align*}
First, the compactness of $H^1(0,T;\bm{H}_{m})$ into
$C([0,T];\bm{H}_{m})$ (recalling that $\bm{H}_{m}$ is a finite-dimensional space) implies that $\Phi^m$ is a compact operator from $C([0,T];\bm{H}_{m})$ into itself. Next, we verify the continuity of $\Phi^m$. Indeed, we see from the estimates \eqref{auphih1}, \eqref{aum} that 
$$ 
\sup_{t\in [0,T]}\|\bm{v}_1^{m}(t)-\bm{v}^{m}_2(t)\| \leq C_T\sup_{t\in [0,T]}\|\bm{u}_1^{m}(t)-\bm{u}_2^{m}(t)\|,
$$
which yields that $\Phi^m$ is a continuous operator from $C([0,T];\bm{H}_{m})$ into itself.

The estimate \eqref{eeenerg} in the proof of Lemma \ref{fp} (see Appendix) yields that
\begin{align}
&\int_0^T \big(\|\nabla \mu^m(t)\|^2 + \|\nabla \sigma^m(t)\|^2\big) dt \notag\\
&\quad \le 
2 M_2 \left(1+M_3T\Big(\sup_{t\in[0,T]}\|\bm{u}^m\|^2+1\Big)e^{M_3T\Big(\sup_{t\in[0,T]}\|\bm{u}^m\|^2+1\Big)}\right),\notag
\end{align}
where the constant $M_2$ is given by \eqref{M3} and depending on the initial data, coefficients of the system, $\Omega$, $m$, $S$, while the constant $M_3$ may depend on coefficients of the system, $\Omega$, $m$. Recalling the estimate \eqref{auvm1}, we then deduce that 
\begin{align}
&\sup_{t\in[0,T]}\|\bm{v}^{m}(t)\|^2\notag\\
&\quad \le M_4 e^{T}\left[\big(\|\bm{v}_{0}\|^2+ 2 M_2\big) +2M_2M_3T\Big(\sup_{t\in[0,T]}\|\bm{u}^m\|^2+1\Big)e^{M_3T\Big(\sup_{t\in[0,T]}\|\bm{u}^m\|^2+1\Big)}\right].\notag
\end{align}
We note that all the constants $M_i, i=1,2,3,4$ are independent of $\sup_{t\in[0,T]}\|\bm{u}^m\|$. Thus, we now choose a sufficiently large constant $\widetilde{M}$ satisfying 
$$\widetilde{M}\geq 4 M_4(\|\bm{v}_{0}\|^2+ 4 M_2).$$
It is easy to check that there exists a sufficiently small $T_m>0$ depending on $\widetilde{M}$ such that 
$$
 M_4 e^{T_m}\left[\|\bm{v}_{0}\|^2+ 2 M_2 +2M_2M_3T_m\big(\widetilde{M}+1\big)e^{M_3T_m\big(\widetilde{M}+1\big)}\right]\leq \widetilde{M}.
 $$
Hence, we can take  
$$\mathbf{K}_m=\Big\{\bm{u}^m\in C([0,T_m];\bm{H}_m): \ \sup_{t\in[0,T_m]}\|\bm{u}^m(t)\|^2\leq \widetilde{M},\quad \bm{u}^{m}(0)= \bm{P}_{\bm{H}_{m}} \bm{v}_{0}\Big\},
$$
which is a closed convex set in $C([0,T_m];\bm{H}_m)$. Then for any $\bm{u}^m\in \mathbf{K}_m$, we see from the above argument that $\bm{v}^m=\Phi^m(\bm{u}^m)\in H^1([0,T_m];\bm{H}_m)\subset\subset C([0,T_m];\bm{H}_m)$ and it satisfies 
\begin{align}
\sup_{t\in[0,T_{m}]}\|\bm{v}^m(t)\|^2 \le\widetilde{M},\notag
\end{align}
namely, $\bm{v}^m\in \mathbf{K}_m$. 
Recall the classical Schauder's fixed point theorem (see e.g., \cite[Chapter 11, Corollary 11.2]{GT}):
\bl\rm    
Let $\mathbf{K}$ be a closed convex set in a Banach space 
$\mathfrak{B}$ and let $\mathcal{T}$ be a continuous mapping of $\mathbf{K}$ into itself such that the image $\mathcal{T}\mathbf{K}$ is precompact. Then $\mathcal{T}$ has a fixed point in $\mathbf{K}$.
\el
\noindent Hence, we are able to conclude that the mapping $\Phi^m$ defined on the space $C([0,T_m];\bm{H}_m)$ has a fixed point $\bm{v}^m$ in the set $\mathbf{K}_m$, and the corresponding $(\varphi^m,\mu^m,\sigma^m)$ are then determined by Lemma \ref{fp}. Besides, uniqueness of the solution $(\bm{v}^m,\varphi^m,\mu^m,\sigma^m)$ to problem (P1) is an easy consequence of the energy method. 

The proof of Proposition \ref{p1} is complete.
$\hfill\blacksquare$

\subsection{Proof of Theorem \ref{main}}
\noindent We are in a position to prove our main result Theorem \ref{main}. 
\subsubsection{Uniform estimates}
\noindent \textbf{First estimate}. From Proposition \ref{p1}, we have
\begin{align}
\|\varphi^m(t)\|_{L^\infty}\leq 1,\quad \forall\,t\in [0,T_m],\label{linfi}
\end{align}
which implies that $\int_\Omega \varPsi(\varphi^m(t))\,dx$ is uniformly bounded due to (H2). Besides, in \eqref{pag4.d}, choosing the test function $\xi=1$, we obtain
\be
\frac{d}{dt}\big(\overline{\varphi}^{m}-c_0\big)+\alpha\big(\overline{\varphi}^{m}-c_0\big)=0,
\notag
\ee
so that 
\be
\overline{\varphi}^{m}(t)=c_0+e^{-\alpha t} \big(\overline{\varphi}_{0}-c_0\big),\quad \forall\, t\in[0,T_m].
\label{mphimean}
\ee
Since $\overline{\varphi}_{0}, c_0\in (-1,1)$ and $\alpha\geq 0$, we have 
\be 
|\overline{\varphi}^{m}(t)|<1,\quad \forall\, t\in [0,T_m].\label{averphi}
\ee
\textbf{Second estimate}. Testing \eqref{ptest3.c} with $\bm{v}^{m}$, \eqref{pag4.d} with $\mu^{m}$, \eqref{pag1.a} with $\partial_t \varphi^{m}$,\ \eqref{pag2.b} with $\sigma^{m}+\chi(1-\varphi^{m})$, adding the resultants together and integrating with respect to time, we obtain 
\begin{align}
& \mathcal{E}^m(t)  + \int_0^t \mathcal{D}^m(\tau)\,d\tau  = \mathcal{E}^m(0) +  \int_0^t \mathcal{R}^m(\tau)\, d\tau,\quad \forall\,t\in[0,T_m].\label{menergy}
\end{align}
where
\begin{align*}
\mathcal{E}^m(t)&=\frac{1}{2}\|\bm{v}^{m}(t)\|^2 
+ \int_\Omega A\varPsi(\varphi^{m}(t))\, dx + \frac{B}{2}\|\nabla \varphi^{m}(t)\|^2
+\frac{1}{2}\|\sigma^{m}(t)\|^2\notag\\
&\quad 
+\int_\Omega \chi\sigma^{m}(t)(1-\varphi^{m}(t))\, dx,\notag\\
\mathcal{D}^m(t)& =\int_\Omega 2\eta(\varphi^{m}(t))|D\bm{v}^{m}(t)|^2\, dx + \|\nabla \mu^{m}(t)\|^2+\|\nabla(\sigma^{m}(t)+\chi(1-\varphi^{m}(t))\|^2,\\
\mathcal{R}^m(t) &=\int_\Omega \big[\sigma^m(t) +\chi(1-\varphi^m(t))\big]\big[-\mathcal{C} h(\varphi^m(t)) \sigma^m(t) +S(x,t)\big]\, dx\\
&\quad -\alpha\int_\Omega(\varphi^m-c_0)\mu^m\,dx.
\end{align*}
The initial energy satisfies
\begin{align}
\mathcal{E}^m(0) &
=\frac{1}{2}\|\bm{P}_{\bm{H}_{m}} \bm{v}_{0}\|^2 
+ \int_\Omega A\varPsi(\varphi_0)\, dx + \frac{B}{2}\|\nabla \varphi_0\|^2
+\frac{1}{2}\|\sigma_0\|^2\notag\\
&\quad 
+\int_\Omega \chi \sigma_0(1-\varphi_0)\, dx\notag\\
&\le C\Big(\|  \bm{v}_{0}\|, \|\varphi_{0}\|_{H^1}, \int_\Omega A\varPsi(\varphi_0)\, dx, \|\sigma_{0}\|, A, B, \chi, \Omega\Big)\notag\\
&:= \mathcal{E}_0,\label{iniE0}
\end{align}
where the constant $\mathcal{E}_0$ is independent of $m$. 
In light of the $L^\infty$ estimate \eqref{linfi}, we deduce that  
\begin{align}
&\int_{\Omega}\frac{1}{2}|\sigma^m|^2+\chi\sigma^m(1-\varphi^m) \ dx\notag\\
&\quad\ge\int_{\Omega}\left(\frac{1}{2}|\sigma^m|^2-2|\chi||\sigma^m|\right)\ dx\notag\\
&\quad\ge \frac{1}{4}\|\sigma^m\|^2-4\chi^2|\Omega|.\label{Lowerbd}
\end{align}
Next, the first term in $\mathcal{R}^m$ can be estimated by 
\begin{align}
&\int_\Omega \big[\sigma^m(t) +\chi(1-\varphi^m(t))\big]\big[-\mathcal{C} h(\varphi^m(t)) \sigma^m(t) +S(x,t)\big]\, dx\notag\\
&\quad \leq (\|\sigma^m\|+|\chi|+|\chi|\|\varphi^m\|)(|\mathcal{C}|\|h(\varphi^m)\|_{L^\infty}\|\sigma^m\|+\|S\|)\notag\\
&\quad \leq C(1+\|\sigma^m\|^2+\|S\|^2),\notag
\end{align}
where the constant $C$ is independent of $m$. On the other hand, using the convexity of $\varPsi_0$ and \eqref{linfi}, we can estimate the second term (see \cite{GGM})
\begin{align}
&-\alpha\int_\Omega(\varphi^m-c_0)\mu^m\,dx\notag\\
&\quad = -\alpha B\|\nabla \varphi^m\|^2-\alpha A\int_\Omega \varPsi'(\varphi^m)(\varphi^m-c_0)\,dx + \alpha\chi\int_\Omega \sigma^m(\varphi^m-c_0)\,dx\notag\\
&\quad \leq -\frac{\alpha B}{2}\|\nabla \varphi^m\|^2
 -\alpha A\int_\Omega \varPsi_0(\varphi^m)\,dx 
 +\alpha A\int_\Omega \varPsi_0(c_0)\,dx \notag\\
&\qquad
+\alpha A\theta_0\int_\Omega\varphi^m(\varphi^m-c_0)\,dx
+ \alpha\chi\int_\Omega \sigma^m(\varphi^m-c_0)\,dx\notag\\
&\quad \leq -\alpha \mathcal{E}^m(t) + \alpha\|\sigma^m\|^2+\frac{\alpha}{2}\|\bm{v}^{m}(t)\|^2 +C,
\label{mueea}
\end{align}
where the constant $C$ may depend on $\alpha$, $A$, $\varPsi_0$, $\theta_0$, $\chi$, $\Omega$, but is independent of $m$. 
 Then we infer from the above estimates and \eqref{menergy} that 
\begin{align}
& \mathcal{E}^m(t)  + \int_0^t \mathcal{D}^m(s)\,ds  \notag\\
&\quad \leq  \mathcal{E}_0 +  C \int_0^t \mathcal{E}^m(s)\, ds + Ct+C\int_0^t\|S(s)\|^2ds,\quad \forall\,t\in[0,T_m],\notag
\end{align}
where the constant $C$ may depend on $\chi$, $A$, $\varPsi$, $h$, $\Omega$, but is independent of $m$. Thus, from Lemma \ref{GronL}, we obtain 
\begin{align}
\mathcal{E}^m(t) + \int_0^t \mathcal{D}^m(s)\,ds \leq Ce^{Ct}\Big(\mathcal{E}_0 + t+\int_0^t\|S(s)\|^2ds\Big), \quad \forall\,t\in[0,T_m].\notag
\end{align} 
As a consequence of the above estimate and \eqref{Lowerbd}, we get the following estimate 
\begin{align}
& \| \bm{v}^{m} (t) \|^2+
 \| \nabla \varphi^{m}(t)\|^2 +  \|  \sigma^{m}(t) \|^2 \notag \\
&\qquad + \int_0^t \left(\|  D\bm{v}^{m}(s) \|^2
+ \| \nabla \mu^{m}(s)\|^2
+ \| \nabla \sigma^{m}(s)\|^2\right) \, ds\notag \\
&\quad \le  C, \quad \forall\,t\in[0,T_m],  
\label{intenergy}
\end{align}
where the constant $C$ depends on $\mathcal{E}_0$, $\Omega$, $t$, $S$, $\eta$, $h$ and coefficients of the system,  but is independent of $m$.  \medskip

\noindent \textbf{Third estimate}. Testing \eqref{pag2.b} by $\xi=1$, we get
\begin{align}
|\overline{\mu}^{m}|&=|\Omega|^{-1}|A(\varPsi'(\varphi^{m}),1)-\chi(\sigma^{m},1)| \notag \\
 & \le |\Omega|^{-1}A\| \varPsi'(\varphi^{m})\|_{L^1}  +|\Omega|^{-\frac12}|\chi| \|  \sigma^{m}\|.
\label{average valuea}
\end{align}
The term $\| \varPsi'(\varphi^{m})\|_{L^1}$ can be estimated as in \cite[Section 3]{GG} (using the argument in \cite{MZ04}) such that 
\begin{equation}
\begin{aligned}
&\| \varPsi'(\varphi^{m})\|_{L^1} \le C \left\|  \nabla \mu^{m}\right \|+C. \nonumber
\end{aligned}
\end{equation}
As a consequence, using Poincar\'{e}'s inequality and \eqref{intenergy}, \eqref{average valuea}, we obtain
\begin{equation}
\left \|  \mu^{m}  \right \|_{L^{2}(0,T_m;H^1(\Omega)) }\le C .\label{mu}
\end{equation}
\textbf{Fourth estimate}. Testing \eqref{pag4.d} with $-\Delta\varphi^{m}$, we get
\begin{align}
&A(\varPsi''(\varphi^{m})\nabla\varphi^{m},\nabla\varphi^{m})+B\left \|  \Delta\varphi^{m}  \right \|^2  \notag \\
&\quad =-A(\varPsi'(\varphi^{m}),\Delta\varphi^{m})+B(\Delta\varphi^{m},\Delta\varphi^{m})\notag \\
&\quad = (\nabla \mu^{m},\nabla\varphi^{m})+(\chi\nabla\sigma^{m},\nabla\varphi^{m}) \notag \\
&\quad\le C(\left \| \nabla \mu^{m}  \right \|\left \| \nabla \varphi^{m}  \right \|+\left \| \nabla \sigma^{m}  \right \|\left \| \nabla \varphi^{m}  \right \|)\notag  \\
&\quad \le C(\left \| \nabla \mu^{m}  \right \|+\left \| \nabla \sigma^{m}  \right \|).
 \label{varP21}
\end{align}
Thus, it follows from (H2) and \eqref{varP21} that
$$B\left \|  \Delta\varphi^{m}  \right \|^2   \le C(\left \| \nabla \mu^{m}  \right \|+\left \| \nabla \sigma^{m}  \right \|+ \left \|  \nabla\varphi^{m}  \right \|^2),$$
which implies 
$$\left \|  \Delta \varphi^{m}\right \|_{L^4(0,T_m;L^2(\Omega))} \le C.$$
By standard elliptic estimates for the Neumann problem, we obtain
\begin{equation}
\left \|  \varphi^{m}\right \|_{L^4(0,T_m;H^2(\Omega))} \le C.
\label{phiH2}
\end{equation}
From Korn's inequality (see e.g., \cite{Ho}), the estimates \eqref{intenergy}, \eqref{mu} and \eqref{phiH2}, we obtain that 
\begin{align}
&\left \|  \bm{v}^{m} \right \|_{L^{\infty}(0,T_m;\bm{L}^2(\Omega))\cap L^{2}(0,T_m;\bm{H}^1(\Omega))} +\left \|  \varphi^{m}\right \|_{L^{\infty}(0,T;H^1(\Omega))\cap L^4(0,T;H^2(\Omega))} \notag \\
& \quad +\left \|  \sigma^{m}\right \|_{L^{\infty}(0,T;L^2(\Omega))\cap L^{2}(0,T;H^1(\Omega) )} +\left \|  \mu^{m}\right \|_{L^{2}(0,T;H^1(\Omega))}\le  C, 
\label{estimate1}
\end{align}
where the constant $C$ is independent of $m$. The above uniform estimate also enables us to extend the local solution $(\bm{v}^m,\varphi^m,\mu^m,\sigma^m)$ from $[0,T_m]$ to the whole interval $[0,T]$. Besides, from \eqref{pag1.a}, \eqref{estimate1}, we have 
\begin{align}
\|\varPsi'(\varphi^m)\|_{L^2(0,T;L^2(\Omega))}\leq C. \label{Psip1}
\end{align}
Consider the following elliptic problem with singular term
\begin{equation}
\label{ELL}
\begin{cases}
-\Delta \varphi^m+\varPsi_0'(\varphi^m)=\mu^m+\theta_0\varphi^m+\chi\sigma^m,\quad &\text{ in }\Omega,\\
\partial_{\bm{n}} \varphi=0, \quad &\text{ on }\partial \Omega.
\end{cases}
\end{equation}
We infer from \cite[Lemma 7.4]{GGW} (see also \cite{A2009,GGM}) that 
\begin{align}
\|\varphi^m\|_{W^{2,q}}+\|\varPsi_0'(\varphi^m)\|_{L^q}
\leq C(1+\|\mu^m\|_{H^1}+\|\varphi^m\|_{H^1}+\|\sigma^m\|_{H^1}),\notag
\end{align}
where $q\geq 2$ if $d=2$ and $q\in[2,6]$ if $d=3$. This fact and \eqref{estimate1} yield that
\begin{align}
\|\varphi^m\|_{L^2(0,T;W^{2,q}(\Omega))}+\|\varPsi'(\varphi^m)\|_{L^2(0,T;L^q(\Omega))}\leq C.
\label{phiW2q}
\end{align}

\medskip

\noindent \textbf{Fifth estimate}. Finally, we derive estimates for the time derivatives.
As $ \varphi^{m}\in  L^4(0,T;H^2(\Omega))$ and $\bm{v}^{m}  \in   L^{\infty}(0,T;\bm{L}^2)$, we infer from the Sobolev embedding theorem that
\begin{equation}
\int_0^T \left \|  \varphi^{m}(t)\bm{v}^{m}(t)\right \|^2\,dt
\leq \|\bm{v}^{m}\|_{L^{\infty}(0,T;\bm{L}^2)}^2 \int_0^T \|  \varphi^{m}(t)\|_{H^2}^2 \,dt  \le C,\nonumber
\end{equation}
and thus in view of \eqref{pag4.d}, it holds
\begin{equation}
\left \|  \partial_{t}\varphi^{m}\right \|^2_{L^{2}(0,T;H^1(\Omega)')} \le C.
\label{phimt}
\end{equation}
Next, by a similar argument as for \cite[(4.15)--(4.18)]{LW}, we obtain, in two dimensions,
\begin{align}
& \left \| \partial_{t}\sigma^{m}\right \|_{ L^2(0,T;H^1(\Omega)')} + \left \| (\bm{v}^{m}  \cdot \nabla)\sigma^{m} \right \|_{ L^2(0,T;H^1(\Omega)')} \le C, \label{sigmt2d}\\
&\left \| \partial_{t}\bm{v}^{m}\right \|_{ L^{2}(0,T;\bm{H}^1_{0,\mathrm{div}}(\Omega)')} + \left \| (\bm{v}^{m}  \cdot \nabla)\bm{v}^{m}\right \|_{ L^{2}(0,T;\bm{H}^1_{0,\mathrm{div}}(\Omega)')} \le C,\label{vmt2d}
\end{align}
and in three dimensions,
\begin{align}
&\left \| \partial_{t}\sigma^{m}\right \|_{ L^{\frac{4}{3}}(0,T;H^1(\Omega)')} + \left \| (\bm{v}^{m}  \cdot \nabla)\sigma^{m} \right \|_{ L^{\frac{4}{3}}(0,T;H^1(\Omega)')} \le C,\label{sigmt3d}\\
&\left \| \partial_{t}\bm{v}^{m} \right \|_{ L^{\frac{4}{3}}(0,T;\bm{H}^1_{0,\mathrm{div}}(\Omega)')} + \left \| (\bm{v}^{m}  \cdot \nabla)\bm{v}^{m}\right \|_{ L^{\frac{4}{3}}(0,T;\bm{H}^1_{0,\mathrm{div}}(\Omega)')} \le C,\label{vmt3d}
\end{align}
where the constant $C$ is independent of $m$.

\subsubsection{Passage to the limit as $m \to +\infty$} \label{limit}
Thanks to the uniform estimates \eqref{estimate1}--\eqref{vmt3d} that are independent of $m$, we are able to apply the compactness argument to conclude that, when letting $m\to +\infty$, there exists a convergent subsequence of the approximate solutions $(\bm{v}^{m},\varphi^{m},\mu^{m},\sigma^{m})$, whose limit denoted by $(\bm{v},\varphi,\mu,\sigma)$ is a global weak solution to problem \eqref{f3.c}--\eqref{ini0}. Since this procedure is standard, we omit the details here. 

The proof of Theorem \ref{main} is complete.
$\hfill\blacksquare$

\section{Uniqueness of Weak Solutions in Dimension Two}\label{uw}
\setcounter{equation}{0}
In this section, we prove the continuous dependence result in Theorem \ref{the2} and the uniqueness of weak solutions in two dimensions (Corollary \ref{uniq}). The main difficulty comes from the variable viscosity, which can be overcome by using the idea in \cite{GMT} for the two dimensional Cahn--Hilliard--Navier--Stokes system. Additional efforts will be made to handle the nutrient equation for $\sigma$ and nonconservation of the mass for $\varphi$. \medskip

\noindent \textbf{Proof of Theorem \ref{the2}}. Let $(\bm{v}_{1},\varphi_{1},\mu_1,\sigma_{1})$ and $(\bm{v}_{2},\varphi_{2},\mu_2,\sigma_{2})$ be two weak solutions to problem \eqref{f3.c}--\eqref{ini0} given by Theorem \ref{main} subject to the initial data $(\bm{v}_{01},\sigma_{01},\varphi_{01})$ and $(\bm{v}_{02},\sigma_{02},\varphi_{02})$, respectively. For simplicity, we assume $S_1=S_2$. Denote the differences
$$
(\bm{v},\varphi,\mu,\sigma)=(\bm{v}_{1}-\bm{v}_{2},\varphi_{1}-\varphi_{2},\mu_1-\mu_2,\sigma_{1}-\sigma_{2}).
$$
By the definition of weak solutions and the divergence free condition, we see that (also recall a similar calculation like \eqref{atest111.a}, \eqref{eqva1}) 
\begin{subequations}
\begin{alignat}{3}
&\left \langle\partial_t  \bm{ v},\bm{\zeta}\right \rangle_{\bm{H}^1_{0,\mathrm{div}}}
-(\bm{ v}_{1} \otimes \bm {v},\nabla\bm{ \zeta})
-(\bm{ v} \otimes \bm {v}_{2},\nabla\bm{ \zeta})
+(2\eta(\varphi_{1}) D\bm{v} ,D\bm{ \zeta}) \notag\\
&\qquad+(  2(\eta(\varphi_{1})-\eta(\varphi_{2})) D\bm{v}_{2} ,D\bm{ \zeta})\notag \\
&\quad=(\nabla \varphi_{1}\otimes\nabla \varphi,\nabla\bm {\zeta})+(\nabla \varphi\otimes\nabla \varphi_{2},\nabla\bm {\zeta}),\label{test33.c} \\
&\left \langle \partial_t \varphi,\xi\right \rangle_{H^1}
-( \varphi \bm{v}_{1},  \nabla\xi)
-( \varphi_{2} \bm{v},  \nabla\xi)\notag\\
&\quad =- (\nabla \mu,\nabla \xi)-\alpha(\varphi, \xi),\label{test11.a} \\
&\  \mu=A\varPsi'(\varphi_{1})-A\varPsi'(\varphi_{2})-B\Delta \varphi-\chi \sigma ,\label{test44.d}\\
&\left \langle\partial_t \sigma,\xi\right \rangle_{H^1}
-(\sigma\bm{v}_{1},  \nabla \xi)
-(\sigma_{2}\bm{v},  \nabla\xi) + (\nabla \sigma, \nabla \xi)\notag\\
&\quad 
=\chi(\nabla \varphi,\nabla \xi)
-\mathcal{C}( h(\varphi_1) \sigma_1- h(\varphi_2)\sigma_2,\xi), \label{test22.b}  
\end{alignat}
\end{subequations}
for all $\bm {\zeta} \in \bm {H}^1_{0,\mathrm{div}}(\Omega)),\ \xi \in H^1(\Omega)$.

First, for the mean value of $\varphi$, we see that (cf. \eqref{mphimean}) 
\begin{align}
\frac{d}{dt}\overline{\varphi}
+\alpha\overline{\varphi}=0, 
\qquad \text{and thus}\ \ \overline{\varphi}(t)=\overline{\varphi}_0e^{-\alpha t}, 
\quad \forall\,t\in[0,T].
\label{diffmean}
\end{align}
From the expression of $\overline{\varphi}$, we also have 
\begin{align}
\frac{1}{2}\frac{d}{dt}\overline{\varphi}^2
+\alpha\overline{\varphi}^2=0 \quad \text{and} \quad \frac{d}{dt}|\overline{\varphi}|
+\alpha|\overline{\varphi}|=0.  
\label{diffmean0}
\end{align}
Taking the test function $\xi=\mathcal{N}(\varphi-\overline{\varphi})$ in \eqref{test11.a}, we obtain
\begin{equation}
\frac{1}{2}\frac{d}{dt}\|\varphi-\overline{\varphi}\|_{V_0'}^2+(\mu,\varphi-\overline{\varphi})+\alpha\|\varphi-\overline{\varphi}\|_{V_0'}^2=I_{1}+I_{2},
\label{diffvmb}
\end{equation}
where
$$I_{1}=(\varphi\bm{v}_{1},\nabla \mathcal{N}(\varphi-\overline{\varphi})),\ \ I_{2}=(\varphi_{2}\bm{v},\nabla \mathcal{N}(\varphi-\overline{\varphi})).$$
From the assumption on $\varPsi$, we have 
\begin{align}
(\mu,\varphi-\overline{\varphi})&=A(\varPsi'(\varphi_{1})-\varPsi'(\varphi_{2}),\varphi)-A(\varPsi'(\varphi_{1})-\varPsi'(\varphi_{2}),\overline{\varphi})\notag\\
&\quad +B(\nabla \varphi,\nabla \varphi)-(\chi \sigma,\varphi-\overline{\varphi}) \notag\\
&\ge B\|\nabla \varphi\|^2-(A|\theta_0-\theta|+\chi^2)\|\varphi\|^2-\frac{1}{2}\|\sigma\|^2-\chi^2|\Omega|\overline{\varphi}^2\notag\\
&\quad -A(\varPsi'(\varphi_{1})-\varPsi'(\varphi_{2}),\overline{\varphi}),\notag
\end{align}
where 
\begin{align}
&(A|\theta_0-\theta|+\chi^2)\|\varphi\|^2\notag\\
&\quad \leq 2 (A|\theta_0-\theta|+\chi^2) \|\varphi-\overline{\varphi}\|^2+ 2(A|\theta_0-\theta|+\chi^2)|\Omega|\overline{\varphi}^2 \notag\\
&\quad =2(A|\theta_0-\theta|+\chi^2)(\nabla \mathcal{N}(\varphi-\overline{\varphi}),\nabla \varphi)+ 2(A|\theta_0-\theta|+\chi^2)|\Omega|\overline{\varphi}^2\notag\\
&\quad \le\frac{B}{2}\|\nabla \varphi\|^2+C\|\varphi-\overline{\varphi}\|_{V_0'}^2+C\overline{\varphi}^2\notag\\
&\quad \le\frac{B}{2}\|\nabla \varphi\|^2+C\|\varphi\|_{(H^1)'}^2.\notag
\end{align}
Thus, from \eqref{diffmean0}, \eqref{diffvmb} and the equivalent norm on $(H^1)'$, we get 
\begin{align}
& \frac{1}{2}\frac{d}{dt}\|\varphi\|_{(H^1)'}^2
+\frac{B}{2}\|\nabla \varphi\|^2+\alpha\overline{\varphi}^2\notag\\
&\le \frac{1}{2}\|\sigma\|^2 +C\|\varphi\|_{(H^1)'}^2
+C\big(\|\varPsi'(\varphi_{1})\|_{L^1}+\|\varPsi'(\varphi_{2})\|_{L^1}\big)|\overline{\varphi}| +I_{1}+I_{2}.
\label{starphi}
\end{align}

Next, taking $\bm{\zeta}=\bm{S}^{-1}\bm{v}$ in \eqref{test33.c}, we get (see \cite[(3.8)]{GMT})
\begin{equation}
\frac{1}{2}\frac{d}{dt}\|\nabla\bm{S}^{-1}\bm{v}\|^2
+(  \eta(\varphi_{1}) D\bm{v} ,\nabla\bm{S}^{-1}\bm{v})=I_{3}+I_{4}+I_{5},
\label{v2}
\end{equation}
where the right-hand side terms are given by 
\begin{align*}
I_{3}&=-((\eta(\varphi_{1})-\eta(\varphi_{2})) D\bm{v}_{2} ,\nabla\bm{S}^{-1}\bm{v}),\\
I_{4}&=(\bm{ v}_{1} \otimes \bm {v},\nabla\bm{S}^{-1}\bm{v})+(\bm{ v} \otimes \bm {v}_{2},\nabla\bm{S}^{-1}\bm{v}),\\
I_{5}&=(\nabla \varphi_{1}\otimes\nabla \varphi,\nabla\bm{S}^{-1}\bm{v})
+(\nabla \varphi\otimes\nabla \varphi_{2},\nabla\bm{S}^{-1}\bm{v}).
\end{align*}
By the property of the Stokes operator, there exists a $p \in L^2(0,T;H^1(\Omega))$, such that
$-\Delta\bm{S}^{-1}\bm{v}+\nabla p =\bm{v}$, a.e. in $\Omega\times(0,T)$ and satisfies (see Lemma \ref{stokes})
$
\|p\|\le C\|\nabla\bm{S}^{-1}\bm{v}\|^{\frac{1}{2}}\|\bm{v} \|^{\frac{1}{2}}$, $\|p\|_{H^1}\le C \| \bm{v}\|$.
Then the following observation was made in \cite[(3.9),(3.11)]{GMT}:
\begin{align}
&(  \eta(\varphi_{1}) D\bm{v},\nabla\bm{S}^{-1}\bm{v})\notag\\
&\quad \ge -(\bm{ v},\eta'(\varphi_{1})D\bm{S}^{-1}\bm{v}\nabla\varphi_{1})+\frac{\eta_{*}}{2}\|\bm{v} \|^2+\frac{1}{2}(\eta'(\varphi_{1})\nabla \varphi_{1}\cdot\bm{v} , p).
\label{OB}
\end{align}
It follows from \eqref{starphi}, \eqref{v2} and \eqref{OB} that 
\begin{align}
&\frac12\frac{d}{dt}\Big(\|\nabla\bm{S}^{-1}\bm{v}\|^2+ \|\varphi\|_{(H^1)'}^2  \Big)+\frac{\eta_{*}}{2}\|\bm{v}\|^2
+\frac{B}{2}\|\nabla \varphi\|^2\notag\\
& \quad \le
\frac{1}{2}\|\sigma\|^2+C\|\varphi\|_{(H^1)'}^2 +C\big(\|\varPsi'(\varphi_{1})\|_{L^1}+\|\varPsi'(\varphi_{2})\|_{L^1}\big)|\overline{\varphi}|+\sum_{j=1}^7 I_j, \label{uniquedif}
\end{align}
where 
$$I_{6}=(\bm{ v},\eta'(\varphi_{1})D\bm{S}^{-1}\bm{v}\nabla\varphi_{1}),\ \ I_{7}=-\frac{1}{2}(\eta'(\varphi_{1})\bm{v}\cdot \nabla \varphi_{1}, p).$$

Taking now the test function $\xi=\mathcal{N}_1\sigma$ in \eqref{test22.b}, we get
\begin{align}
&\frac12\frac{d}{dt}\|\sigma\|_{(H^1)'}^2
 + \|\sigma\|^2=\sum_{j=8}^{14}I_j,\label{uniquedifs}
\end{align}
where  
\begin{align*}
I_8&=(\sigma\bm{v}_{1},  \nabla \mathcal{N}_1\sigma),\qquad\quad  
I_9=(\sigma_{2}\bm{v},  \nabla\mathcal{N}_1\sigma),\\
I_{10}&=(\sigma, \mathcal{N}_1\sigma),\quad \qquad \qquad 
I_{11}=\chi(\nabla \varphi,\nabla \mathcal{N}_1\sigma),\\
I_{12}&=-\mathcal{C}( h(\varphi_1) \sigma,\mathcal{N}_1\sigma),\quad \, 
I_{13}=-\mathcal{C}( (h(\varphi_1)- h(\varphi_2))\sigma_2,\mathcal{N}_1\sigma).
\end{align*}
It remains to estimate the reminder terms $I_1,...,I_{13}$ on the  right-hand side. Denote 
 \begin{align}
 W(t):=\|\nabla\bm{S}^{-1}\bm{v}(t)\|^2+\|\varphi(t)\|_{(H^1)'}^2+\|\sigma(t)\|_{(H^1)'}^2+|\overline{\varphi}(t)|.
 \label{W}
 \end{align}
The terms $I_1,I_2,I_4,I_5,I_6,I_7$ in \eqref{uniquedif} can be estimated as in \cite[Section 3]{GMT} with minor modifications such that 
\begin{align*}
I_{1}
&= ((\varphi-\overline{\varphi})\bm{v}_{1},\nabla \mathcal{N}(\varphi-\overline{\varphi}))\notag\\
&\le\frac{B}{20}\|\nabla\varphi\|^2+C\|\bm{v}_{1}\|_{\bm{L}^3}^2\|\varphi-\overline{\varphi}\|_{V_0'}^2\\
&\le \frac{B}{20}\|\nabla\varphi\|^2+C\|\nabla \bm{v}_{1}\|^2\|\varphi\|_{(H^1)'}^2,\\
I_{2}
&\le\frac{\eta_*}{20}\|\bm{v}\|^2+C\|\varphi-\overline{\varphi}\|_{V_0'}^2\le \frac{\eta_*}{20}\|\bm{v}\|^2+C\|\varphi\|_{(H^1)'}^2,\\
I_{4}
&\le\frac{\eta_*}{20}\|\bm{v}\|^2
+C(\|\nabla \bm{v}_{1}\|^2+
\|\nabla \bm{v}_{2}\|^2)\|\nabla\bm{S}^{-1}\bm{v}\|^2,\\
I_{5}
&\le\frac{B}{20}\|\nabla\varphi\|^2+
C(\|\nabla \varphi_{1}\|_{\bm{L}^{\infty}}^2+\|\nabla \varphi_{2}\|_{\bm{L}^{\infty}}^2)\|\nabla\bm{S}^{-1}\bm{v}\|^2\\
&\le\frac{B}{20}\|\nabla\varphi\|^2+
C(\|\varphi_{1}\|_{W^{2,3}}^2+\|\varphi_{2}\|_{W^{2,3}}^2)\|\nabla\bm{S}^{-1}\bm{v}\|^2,\\
I_{6}&\le\frac{\eta_*}{20}\|\bm{v}\|^2+C\|\nabla \varphi_{1}\|_{\bm{L}^{\infty}}^2\|\nabla\bm{S}^{-1}\bm{v}\|^2\\
&\le\frac{\eta_*}{20}\|\bm{v}\|^2+C\|\varphi_{1}\|_{W^{2,3}}^2\|\nabla\bm{S}^{-1}\bm{v}\|^2,\\
I_{7}&\le\frac{\eta_*}{20}\|\bm{v}\|^2+C\|\varphi_{1}\|_{H^{2}}^4\|\nabla\bm{S}^{-1}\bm{v}\|^2.
\end{align*}
For $I_3$, we recall the following result (see \cite[Proposition C.2]{GMT}):
\begin{lemma}\label{log}\rm 
Let $\Omega$ be a bounded domain in $\mathbb{R}^2$ with smooth boundary. Assume
that $f, h \in H^1(\Omega)$ and $\bm{g}\in \bm{H}^1(\Omega)$. Then, there exists a positive constant $C$ such
that 
\be
\|f\bm{g}\| \le C \|f\|_{H^1}(\|\bm{g}\|+\|h\|) \left[\ln\Big(e\frac{\|\bm{g}\|_{\bm{H}^1}+\|h\|_{H^1}}{\|\bm{g}\|+\|h\|}\Big)\right]^{\frac{1}{2}}.\notag
\ee
\end{lemma}
\noindent Then using the assumption (H1) and taking 
$$f=\varphi-\overline{\varphi},\quad  h=\mathcal{N}^\frac12(\varphi-\overline{\varphi})+\mathcal{N}_1^\frac12\sigma+|\overline{\varphi}|^\frac12,\quad  \bm{g}=\nabla \bm{S}^{-1}\bm{v}$$ in Lemma \ref{log}, we can deduce that 
\begin{align}
I_{3}
&\le \eta_0\|D\bm{v}_{2}\|\|(\varphi-\overline{\varphi})\nabla \bm{S}^{-1}\bm{v}\|+ \eta_0\|D\bm{v}_{2}\|\|\overline{\varphi}\nabla \bm{S}^{-1}\bm{v}\| \notag\\
&\le C\|D\bm{v}_{2}\|\|\varphi-\overline{\varphi}\|_{H^1}\big(\|\nabla\bm{S}^{-1}\bm{v}\|+\|\varphi-\overline{\varphi}\|_{V_0'}+\|\sigma\|_{(H^1)'}+|\Omega|^\frac12|\overline{\varphi}|^\frac12\big)\notag\\
&\quad \times 
\left[\ln\left(e\frac{\|\nabla\bm{S}^{-1}\bm{v}\|_{\bm{H}^1}+\|\mathcal{N}^\frac12(\varphi-\overline{\varphi})\|_{H^1}+\|\mathcal{N}_1^\frac12\sigma\|_{H^1}+|\Omega|^\frac12|\overline{\varphi}|^\frac12}{\|\nabla\bm{S}^{-1}\bm{v}\|+\|\varphi-\overline{\varphi}\|_{V_0'} +\|\sigma\|_{(H^1)'}+|\Omega|^\frac12|\overline{\varphi}|^\frac12}\right)\right]^\frac12\notag\\
&\quad + \eta_0|\overline{\varphi}|\|D\bm{v}_{2}\|\nabla \bm{S}^{-1}\bm{v}\|\notag\\
&\le C\|D\bm{v}_{2}\|\|\nabla \varphi\|\big(\|\nabla\bm{S}^{-1}\bm{v}\|+\|\varphi-\overline{\varphi}\|_{V_0'}+\|\sigma\|_{(H^1)'}+|\Omega|^\frac12|\overline{\varphi}|^\frac12\big)\notag\\
&\quad \times 
\left(\ln\frac{Ce(\|\bm{v}\|+\|\varphi\|+\|\sigma\|+1)}{\|\nabla\bm{S}^{-1}\bm{v}\|+\|\varphi\|_{(H^1)'} +\|\sigma\|_{(H^1)'}+ |\overline{\varphi}|^\frac12}\right)^\frac12
+ \eta_0|\overline{\varphi}|\|D\bm{v}_{2}\|\|\nabla \bm{S}^{-1}\bm{v}\|\nonumber\\
&\le \frac{B}{20}\|\nabla\varphi\|^2+ C\|D\bm{v}_{2}\|^2W(t)\ln \left(\hat{C}e\frac{\|\bm{v}\|^2+\|\varphi\|^2+\|\sigma\|^2+1}{W(t)}\right)+|\overline{\varphi}|^2,\notag
\end{align}
where we used the fact $e^{-\frac{\alpha}{2}t}\geq e^{-\alpha t}$ for $t\in[0,T]$. 
We remark that from its definition, the logarithmic term satisfies
\begin{align}
\ln \left(\hat{C}e\frac{\|\bm{v}\|^2+\|\varphi\|^2+\|\sigma\|^2+1}{W(t)}\right)\geq 1, \notag
\end{align}
provided that we choose the constant $\hat{C}$ properly large. On the other hand, under this choice, and thanks to the boundedness of $\|\bm{v}\|_{L^\infty(0,T;\bm{L}^2(\Omega))}$, $\|\varphi\|_{L^\infty(0,T;L^2(\Omega))}$, $\|\sigma\|_{L^\infty(0,T;L^2(\Omega))}$ (cf. Theorem \ref{main}), it also holds 
\begin{align}
1\le \ln \left(\hat{C}e\frac{\|\bm{v}\|^2+\|\varphi\|^2+\|\sigma\|^2+1}{W(t)}\right)\leq \ln\left(\frac{\widetilde{C}}{W(t)}\right)=-\ln\left(\frac{W(t)}{\widetilde{C}}\right), \label{CCC}
\end{align}
for some constant $\widetilde{C}>0$. 

Next, we estimate the terms $I_8,...,I_{13}$ in the $\sigma$-equation.
\begin{align}
I_8
&\le C\|\bm{v}_{1}\|_{\bm{L}^4}\|\sigma\|\|\nabla \mathcal{N}_1\sigma\|_{\bm{L}^4}\notag\\
&\le C\|\bm{v}_{1}\|^{\frac{1}{2}}\|\nabla\bm{v}_{1}\|^{\frac{1}{2}}\|\sigma\|
\|\nabla \mathcal{N}_1\sigma\|^{\frac{1}{2}}\|\nabla \mathcal{N}_1\sigma\|_{\bm{H}^1}^{\frac{1}{2}}\notag\\
&\le \frac{1}{8}\|\sigma\|^2 +C\|\nabla\bm{v}_{1}\|^2\|\sigma\|_{(H^1)'}^2,\notag
\end{align}
\begin{align}
I_9
&\le C\|\sigma_{2}\|_{L^4}\|\bm{v}\|\|\nabla \mathcal{N}_1\sigma\|_{\bm{L}^4}\notag\\
&\le C\|\sigma_2\|^\frac12\|\sigma_2\|_{H^1}^\frac12\|\bm{v}\|
\|\nabla \mathcal{N}_1\sigma\|^\frac12
\|\nabla \mathcal{N}_1\sigma\|_{\bm{H}^1}^\frac12  \notag\\
&\le \frac{\eta_*}{20} \|\bm{v}\|^2
+\frac{1}{8}\|\sigma\|_{L^2}^2
+C\|\sigma_{2}\|_{H^1}^2 \|\sigma\|_{(H^1)'}^2,\notag
\end{align}
\begin{align}
&I_{10}+I_{11}+I_{12}\notag\\
&\quad \le \|\sigma\|_{(H^1)'}^2+|\chi|\|\nabla \varphi\|\|\nabla \mathcal{N}_1\sigma\|+|\mathcal{C}|\|h(\varphi_1)\|_{L^\infty}\|\sigma\|_{(H^1)'}^2\notag\\
&\quad \le \frac{B}{20}\|\nabla\varphi\|^2
+ C\|\sigma\|_{(H^1)'}^2,\notag
\end{align}
and 
\begin{align}
I_{13} 
&=-\mathcal{C}\left(\int_0^1 h'(s\varphi_1+(1-s)\varphi_2)\,ds\,\varphi\sigma_2,\mathcal{N}_1\sigma\right)\notag\\ 
&\le C\int_0^1 \|h'(s\varphi_1+(1-s)\varphi_2)\|_{L^\infty}\,ds\|\varphi\|_{L^4}\|\sigma_2\|\|\mathcal{N}_1\sigma\|_{L^4}\notag\\
&\le C\|\varphi-\overline{\varphi}\|_{V_0'}^\frac14
\|\varphi-\overline{\varphi}\|_{H^1}^\frac34\|\mathcal{N}_1\sigma\|_{H^1}+ C\|\overline{\varphi}\|_{L^4}\|\mathcal{N}_1\sigma\|_{H^1}\notag\\
&\le \frac{B}{20}\|\nabla\varphi\|^2 +C\|\varphi-\overline{\varphi}\|_{V_0'}^2+C|\overline{\varphi}|^2
+C\|\sigma\|_{(H^1)'}^2\notag\\
&\leq \frac{B}{20}\|\nabla\varphi\|^2 +C\|\varphi\|_{(H^1)'}^2
+C\|\sigma\|_{(H^1)'}^2.\notag
\end{align}
Collecting the above estimates, we deduce from \eqref{diffmean0}, \eqref{uniquedif}, \eqref{uniquedifs} and \eqref{W} that 
\begin{align}
\frac{d}{dt}W(t)+ \frac{\eta_*}{2} \|\bm{v}\|^2 + \frac{B}{2}\|\nabla\varphi\|^2 + \|\sigma\|^2
\leq -C_1Z(t)W(t)\ln\left(\frac{W(t)}{\widetilde{C}}\right),
\label{uniA}
\end{align}
where $W(t)$ is defined in \eqref{W}, 
\begin{align}
Z(t)&=\|\nabla \bm{v}_{1}(t)\|^2+\|\nabla \bm{v}_{2}(t)\|^2+\|\varphi_{1}(t)\|_{W^{2,3}}^2+\|\varphi_{2}(t)\|_{W^{2,3}}^2\notag\\
&\quad +\|\varphi_{1}(t)\|_{H^{2}}^4+\|\varPsi'(\varphi_{1})\|_{L^1}+\|\varPsi'(\varphi_{2})\|_{L^1}+\|\sigma_2(t)\|_{H^1}^2+1.\label{Z}
\end{align}
and $C_1$, $\widetilde{C}$ are constants depending on the initial data, $\Omega$, and the coefficients of the system. 

Recalling that $Z(t)\in L^1(0,T)$ for any $T>0$ (see Theorem \ref{main} and \eqref{Psip1}), we conclude from \eqref{uniA} that 
\be
W(t)\le \widetilde{C}\left(\frac{W(0)}{\widetilde{C}}\right)^{\exp\left(-C_1\int_{0}^{t}Z(s)\, ds\right)},\ \ \forall\, t\in[0,T].
\label{contiuniq}
\ee
The uniqueness of weak solutions to problem \eqref{f3.c}--\eqref{ini0} is an immediate consequence of the continuous dependence estimate \eqref{contiuniq}.

The proof is complete. $\hfill\blacksquare$ 

\begin{remark}
Theorem \ref{the2} extends the previous results in \cite{MT,GG,GMT} to a more general context.
The proof of Theorem \ref{the2} also enables us to obtain the uniqueness of global weak solutions to problem \eqref{f3.c}--\eqref{ini0} with unmatched viscosities and a regular polynomial potential like in \cite{LW} (cf. \cite[Remark 3.3]{GMT} for further details).
\end{remark}

\section{Appendix}
\setcounter{equation}{0}
\noindent In the Appendix, we sketch the proof of Lemma \ref{fp} for the well-posedness of the auxiliary problem \eqref{g1.a}--\eqref{g6.h}.\medskip 

\noindent \textbf{Step 1. The regularized problem}.
Concerning the singular potential $\varPsi$ satisfying (H2), without loss of generality, we assume that  $\varPsi_0(0)=0$. Then we may approximate the singular part $\varPsi_0'$, e.g., as in \cite{MT}:
\begin{equation}
\varPsi'_ {0,\epsilon}(r)=\left\{
\begin{aligned}
&\varPsi_0'(-1+\epsilon) + \varPsi_0''(-1+\epsilon)(r+1-\epsilon),\quad r<-1+\epsilon,\\
&\varPsi_0'(r),\qquad\qquad\qquad\qquad\qquad\qquad\qquad|r|\leq 1-\epsilon,\\
&\varPsi_0'(1-\epsilon) + \varPsi_0''(1-\epsilon)(r-1+\epsilon),\qquad\ \  r>1-\epsilon,
\end{aligned}
\right.\label{vPsi}
\end{equation}
for sufficiently small $\epsilon\in(0,\epsilon_0)$ (recall assumption (H2)).  
Define 
$$
\varPsi_{0,\epsilon}(r)=\int_0^{r}\varPsi'_{0,\epsilon}(s) \, ds,\quad \varPsi_\epsilon(r)=\varPsi_{0,\epsilon}(r)-\frac{\theta_0}{2}r^2.
$$
We can verify that $\varPsi_ {0,\epsilon}''(r)\ge \theta>0$ and $\varPsi_ {0,\epsilon}(r)\geq -L$ for $r\in \mathbb{R}$, where $L>0$ is a constant independent of $\epsilon$. Moreover, it holds $\varPsi_{0,\epsilon}(r)\leq \varPsi_{0}(r)$ for $r\in [-1,1]$ (see e.g., \cite{FG12}).

We now introduce the following regularized problem of \eqref{g1.a}--\eqref{g6.h}:
\begin{subequations}
	\begin{alignat}{3}
	&\left \langle\partial_t  \varphi_{\epsilon}^{m},\xi\right \rangle_{{H}^1}+((\bm{u}^{m} \cdot \nabla) \varphi_{\epsilon}^{m},\xi)\notag\\
	&\quad =-(\nabla \mu_{\epsilon}^{m},\nabla\xi)-\alpha(\varphi^m_\epsilon-c_0,\xi),\qquad \qquad\qquad\quad\qquad \, \textrm{a.e.\ in}\ (0,T),\label{g1.ar} \\
	&\mu_{\epsilon}^{m}=A\varPsi'_{\epsilon}(\varphi_{\epsilon}^{m})-B\Delta \varphi_{\epsilon}^{m}-\chi \sigma_{\epsilon}^{m}, \qquad \qquad \qquad \qquad \qquad\rm{a.e.\ in}\ \Omega\times(0,T),\label{g4.dr} \\
	&\left \langle\partial_t  \sigma_{\epsilon}^{m},\xi\right \rangle_{{H}^1}+((\bm{u}^{m} \cdot \nabla) \sigma_{\epsilon}^{m},\xi)+(\nabla \sigma_{\epsilon}^m,\nabla \xi)\notag\\
&\quad = \chi ( \nabla \varphi_{\epsilon}^m,\nabla \xi)-(\mathcal{C} h(\varphi_{\epsilon}^m) \sigma_{\epsilon}^m,\xi) + (S,\xi),\quad\, \qquad \qquad \text{a.e. in}\ (0,T), \label{g2.br}\\
	&\varphi_{\epsilon}^{m}(0)=\varphi_{0},\quad  \sigma_{\epsilon}^{m}(0)=\sigma_{0},\ \ \qquad \qquad \qquad \qquad \qquad \qquad\textrm{in}\ \Omega,  \label{g6.hr}
	\end{alignat}
\end{subequations}
for all $\xi \in H^1(\Omega)$.\medskip

\noindent  \textbf{Step 2. Uniform estimates}. 
The existence and uniqueness of global weak solutions on $[0,T]$ to the regularized problem \eqref{g1.ar}--\eqref{g6.hr} can be proved by using a suitable Galerkin method similar to that in \cite{GL17e}. Below we only derive some uniform estimates with respect to the parameter $\epsilon$ and omit the other details. The process can be made rigorous by the Galerkin scheme.\medskip
  
\noindent \textbf{First estimate}. 
Testing \eqref{g1.ar} by $1$, we get
\be
\frac{d}{dt}(\overline{\varphi}^{m}_\epsilon-c_0)+\alpha(\overline{\varphi}_\epsilon^m-c_0)=0,\notag
\ee
so that 
\be
\overline{\varphi}^{m}_\epsilon(t)=c_0+e^{-\alpha t}(\overline{\varphi}_{0}-c_0),\quad \forall\, t\in [0,T].
\label{averphie}
\ee
\noindent \textbf{Second estimate}. 
Testing \eqref{g1.ar} by $\varphi_{\epsilon}^{m}$ and \eqref{g2.br} by $\sigma^m_\epsilon$, adding the resultants together, we get
\begin{align}
&\frac12 \frac{d}{dt}\Big(\| \varphi_{\epsilon}^{m}\|^2+\|\sigma_{\epsilon}^{m}\|^2\Big)
+B\|\Delta\varphi_{\epsilon}^{m}\|^2+\|\nabla \sigma_{\epsilon}^{m}\|^2 + \alpha \|\varphi_\epsilon^m\|^2 \notag \\
&\quad=\int_{\Omega}\big(A\varPsi'_{\epsilon}(\varphi_{\epsilon}^{m})\Delta\varphi_{\epsilon}^{m} -2\chi\Delta\varphi_{\epsilon}^{m}\sigma^m_\epsilon  -\mathcal{C} h(\varphi_{\epsilon}^m) |\sigma_{\epsilon}^m|^2 + S\sigma_{\epsilon}^m +\alpha c_0 \varphi_{\epsilon}^m\big)\, dx .\label{energy e}
\end{align}
The first term on the right-hand side of \eqref{energy e} can be estimated as follows  
\begin{align}
&\int_{\Omega}A\varPsi'_{\epsilon}(\varphi_{\epsilon}^{m})\Delta\varphi_{\epsilon}^{m} \ dx\notag\\
&\quad=A\int_{\Omega}(\varPsi'_{0,\epsilon}(\varphi_{\epsilon}^{m})-\theta_{0}\varphi_{\epsilon}^{m})\Delta\varphi_{\epsilon}^{m} \, dx\notag\\
&\quad=- A\int_{\Omega}\varPsi_{0,\epsilon}''(\varphi_{\epsilon}^{m})|\nabla\varphi_{\epsilon}^{m}|^2 \, dx
-A\int_{\Omega}\theta_{0}\varphi_{\epsilon}^{m}\Delta\varphi_{\epsilon}^{m} \ dx\notag\\
&\quad\le \frac{B}{4}\|\Delta\varphi_{\epsilon}^{m}\|^2  +\frac{A^2\theta_0^2}{B}\|\varphi_{\epsilon}^{m}\|^2.
\notag
\end{align}
Next, using (H3) and Young's inequality, we get
\begin{align}
&\int_{\Omega}-2\chi\Delta\varphi_{\epsilon}^{m}\sigma_{\epsilon}^{m}dx \le \frac{B}{4}\|\Delta\varphi_{\epsilon}^{m}\|^2 +\frac{4\chi^2}{B}\|\sigma_{\epsilon}^{m}\|^2,\notag\\
&\int_\Omega\big(-\mathcal{C} h(\varphi_{\epsilon}^m) |\sigma_{\epsilon}^m|^2 + S\sigma_{\epsilon}^m +\alpha c_0\varphi_{\epsilon}^m\big)\, dx
\leq C\|\sigma_{\epsilon}^{m}\|^2+\frac12\|S\|^2+C.\notag
\end{align}
Hence, from \eqref{energy e} we see that
\begin{align}
& \frac{d}{dt}\Big(\| \varphi_{\epsilon}^{m}\|^2+\|\sigma_{\epsilon}^{m}\|^2\Big)
+B\|\varphi_{\epsilon}^{m}\|_{H^2}^2+2\|\sigma_{\epsilon}^{m}\|^2_{H^1} \notag \\
&\quad\le C\big(\| \varphi_{\epsilon}^{m}\|^2+ \|\sigma_{\epsilon}^{m}\|^2\big) +\|S\|^2+C.\label{espsL2}
\end{align}

\noindent \textbf{Third estimate}. Next, testing \eqref{g1.ar} with $\mu_{\epsilon}^{m}$,  \eqref{g4.dr} with $\partial_{t}\varphi_{\epsilon}^{m}$, we obtain
\begin{align}
&\frac{d}{dt} \int_{\Omega}  \Big(A\varPsi_{\epsilon}(\varphi_{\epsilon}^{m})
+\frac{B}{2}|\nabla \varphi_{\epsilon}^{m}|^2\Big) \, dx 
+\|\nabla \mu_{\epsilon}^{m}\|^2 \notag \\
&\quad=\int_{\Omega}\big[\chi\sigma_{\epsilon}^{m}\partial_{t}\varphi_{\epsilon}^{m} -(\bm{u}^{m}\cdot \nabla)\varphi_{\epsilon}^{m}\mu_{\epsilon}^{m}-\alpha(\varphi^m_\epsilon -c_0) \mu^m_\epsilon  \big]\, dx.
\label{energyA}
\end{align}
Since $\varphi_{\epsilon}^{m}$ satisfies equation \eqref{g1.ar}, we have 
\begin{align}
\|\partial_{t}\varphi_{\epsilon}^{m}\|_{(H^1)'}
&\le\|\bm{u}^{m}\|_{\bm{L}^3} \|\nabla\varphi_{\epsilon}^{m}\|
+ \|\nabla\mu _{\epsilon}^{m}\|\notag\\
&\le C\|\bm{u}^{m}\| \|\nabla\varphi_{\epsilon}^{m}\|
+\|\nabla\mu _{\epsilon}^{m}\|,\nonumber
\end{align}
which implies 
\begin{align}
\int_{\Omega}\chi\sigma_{\epsilon}^{m}\partial_{t}\varphi_{\epsilon}^{m}\ dx
&\le |\chi|\|\partial_{t}\varphi_{\epsilon}^{m}\|_{(H^1)'}\|\sigma^m_\epsilon\|_{H^1}\notag\\
&\le \frac{1}{4}\|\nabla \mu_{\epsilon}^{m}\|^2 
+ C\chi^2\|\bm{u}^{m}\|^2 \|\nabla\varphi_{\epsilon}^{m}\|^2
+ (1+\chi^2)\|\sigma^m_\epsilon\|_{H^1}^2.\notag
\end{align}
Besides, using Poincar\'{e}'s inequality and \eqref{averphie}, we get 
\begin{align}
-\int_{\Omega}(\bm{u}^{m}\cdot \nabla)\varphi_{\epsilon}^{m}\mu_{\epsilon}^{m}\, dx
&=\int_{\Omega}(\bm{u}^{m}\cdot \nabla)\mu_{\epsilon}^{m}\varphi_{\epsilon}^{m}\, dx \notag\\
&\le \frac{1}{4}\|\nabla \mu_{\epsilon}^{m}\|^2 +
C\|\bm{u}^m\|_{\bm{L}^3}^2\|\varphi_{\epsilon}^{m}\|_{L^6}^2 \notag\\
&\le \frac{1}{4}\|\nabla \mu_{\epsilon}^{m}\|^2 + C\|\bm{u}^m\|^2\Big(\|\nabla \varphi_{\epsilon}^{m}\|^2+1\Big).\notag
\end{align}
Then similar to \eqref{mueeb}, we obtain 
\begin{align}
&-\alpha\int_\Omega(\varphi^m_{\epsilon}-c_0)\mu^m_{\epsilon}\,dx\notag\\
&\quad = -\alpha B\|\nabla \varphi_{\epsilon}^m\|^2-\alpha A\int_\Omega \varPsi'(\varphi_{\epsilon}^m)(\varphi_{\epsilon}^m-c_0)\,dx + \alpha\chi\int_\Omega \sigma_{\epsilon}^m(\varphi_{\epsilon}^m-c_0)\,dx\notag\\
&\quad \leq -\frac{\alpha B}{2}\|\nabla \varphi_{\epsilon}^m\|^2
 -\alpha A\int_\Omega \varPsi_{0,\epsilon}(\varphi_{\epsilon}^m)\,dx 
 +\alpha A\int_\Omega \varPsi_{0,_{\epsilon}}(c_0)\,dx \notag\\
&\qquad
+\alpha A\theta_0\int_\Omega\varphi_{\epsilon}^m(\varphi_{\epsilon}^m-c_0)\,dx
+ \alpha\chi\int_\Omega \sigma_{\epsilon}^m(\varphi_{\epsilon}^m-c_0)\,dx\notag\\
&\quad \leq -\alpha \left(\frac{B}{2}|\nabla \varphi_{\epsilon}^{m}|^2+A\int_\Omega \varPsi_{\epsilon}(\varphi_{\epsilon}^{m})\,dx\right)+ C(\|\varphi_{\epsilon}^m\|^2+ \|\sigma_{\epsilon}^m\|^2)+ C\notag\\
&\quad \leq C(\|\varphi_{\epsilon}^m\|^2+ \|\sigma_{\epsilon}^m\|^2)+ C,
\label{mueeb}
\end{align}
where we use the fact
$$\int_\Omega \varPsi_\epsilon(\varphi^m_\epsilon)\,dx\geq -L|\Omega|-\frac{|\theta_0|}{2}\|\varphi^m_\epsilon\|^2.$$
From the above estimates, we infer from \eqref{energyA} that 
\begin{align}
&\frac{d}{dt} \int_{\Omega}  \Big(A\varPsi_{\epsilon}(\varphi_{\epsilon}^{m})+\frac{B}{2}|\nabla \varphi_{\epsilon}^{m}|^2\Big) \, dx + \frac12 \|\nabla \mu_{\epsilon}^{m}\|^2\notag \\
&\quad\leq  C(1+\chi^2)\|\bm{u}^m\|^2\Big(\|\nabla \varphi_{\epsilon}^{m}\|^2+1\Big) +(1+\chi^2)\|\sigma^m_\epsilon\|_{H^1}^2\notag\\
&\qquad 
+C(\|\varphi_{\epsilon}^m\|^2+ \|\sigma_{\epsilon}^m\|^2)+ C.
\label{energyA1}
\end{align}
\textbf{Fourth estimate}. Multiplying \eqref{espsL2} by $(1+\chi^2)(1+A|\theta_0|)$ and adding the result with \eqref{energyA1}, we get
\begin{align}
&\frac{d}{dt} \int_{\Omega}  \Big[A\varPsi_{\epsilon}(\varphi_{\epsilon}^{m})+\frac{B}{2}|\nabla \varphi_{\epsilon}^{m}|^2+(1+\chi^2)(1+A|\theta_0|)\big(| \varphi_{\epsilon}^{m}|^2+|\sigma_{\epsilon}^{m}|^2\big)\Big] \, dx \notag\\
&\qquad + \frac12 \|\nabla \mu_{\epsilon}^{m}\|^2 
+(1+\chi^2)B\|\varphi_{\epsilon}^{m}\|_{H^2}^2+(1+\chi^2)\|\sigma_{\epsilon}^{m}\|^2_{H^1} \notag \\
&\quad\leq  C(1+\chi^2)\|\bm{u}^m\|^2\Big(\|\nabla \varphi_{\epsilon}^{m}\|^2+1\Big)+ C\big(\| \varphi_{\epsilon}^{m}\|^2+ \|\sigma_{\epsilon}^{m}\|^2+1\big)\nonumber\\
&\qquad  +(1+\chi^2)(1+A|\theta_0|)\|S\|^2.
\label{energyA2}
\end{align}
Hence,  we deduce from \eqref{energyA2} that 
 \begin{align}
&\frac{d}{dt} \widehat{\mathcal{E}}^m_\epsilon(t)  + \frac12 \|\nabla \mu_{\epsilon}^{m}\|^2 
+B\|\varphi_{\epsilon}^{m}\|_{H^2}^2+\|\sigma_{\epsilon}^{m}\|^2_{H^1} \notag \\
&\quad\leq  C\big(\|\bm{u}^m\|^2+1) \widehat{\mathcal{E}}^m_\epsilon(t)  + (1+\chi^2)(1+A|\theta_0|)\|S\|^2,
\label{energyA3}
\end{align}
where 
\begin{align}
 \widehat{\mathcal{E}}^m_\epsilon(t)& =\int_{\Omega}  \Big[A\varPsi_{\epsilon}(\varphi_{\epsilon}^{m})+\frac{B}{2}|\nabla \varphi_{\epsilon}^{m}|^2+(1+\chi^2)(1+A|\theta_0|)\big(| \varphi_{\epsilon}^{m}|^2+|\sigma_{\epsilon}^{m}|^2\big)\Big] \, dx+AL|\Omega|\notag\\
 &\geq \int_{\Omega}  \Big(\frac{B}{2}|\nabla \varphi_{\epsilon}^{m}|^2+| \varphi_{\epsilon}^{m}|^2+|\sigma_{\epsilon}^{m}|^2\Big) \, dx.
 \label{ehat}
\end{align}
Besides, the initial datum satisfies 
\begin{align}
\widehat{\mathcal{E}}^m_\epsilon(0)
&=\int_{\Omega}  \Big[A\varPsi_{\epsilon}(\varphi_0)+\frac{B}{2}|\nabla \varphi_0|^2+(1+\chi^2)(1+A|\theta_0|)\big(| \varphi_0|^2+|\sigma_0|^2\big)\Big] \, dx+AL|\Omega|\notag\\
&\leq \int_{\Omega}  \Big[A\varPsi(\varphi_0)+\frac{B}{2}|\nabla \varphi_0|^2+(1+\chi^2)(1+A|\theta_0|)\big(| \varphi_0|^2+|\sigma_0|^2\big)\Big] \, dx+AL|\Omega|\notag\\
&:=M_1,\label{iniM2}
\end{align}
where the positive constant $M_1$ depends on the initial data, coefficients of the problem, $\Omega$, but is independent of the parameter $\epsilon$. 

It follows from \eqref{energyA3}, Lemma \ref{GronL} and \eqref{iniM2} that 
\begin{align}
&\widehat{\mathcal{E}}^m_\epsilon(t)  
+ \int_0^t \left(\frac12 \|\nabla \mu_{\epsilon}^{m}(s)\|^2 
+ B\|\varphi_{\epsilon}^{m}(s)\|_{H^2}^2
+ \|\sigma_{\epsilon}^{m}(s)\|^2_{H^1} \right) ds\notag\\
&\quad\leq M_2 \left(1+Ct\Big(\sup_{t\in[0,T]}\|\bm{u}^m\|^2+1\Big)e^{Ct\left(\sup_{t\in[0,T]}\|\bm{u}^m\|^2+1\right)}\right),
\label{eeenerg}
\end{align}
for all $t\in [0,T]$, where
\begin{align}
M_2=M_1+(1+\chi^2)(1+A|\theta_0|)\|S\|_{L^2(0,T;L^2(\Omega))}^2.
\label{M3}
\end{align}
In light of \eqref{ehat}, \eqref{eeenerg}, we obtain
\begin{align}
&\left \| \varphi_{\epsilon}^{m}\right \|_{L^{\infty}(0,T;H^1(\Omega))\cap L^2(0,T;H^2_{N}(\Omega))}
+\left \| \sigma_{\epsilon}^{m}\right \|_{L^{\infty}(0,T;L^2(\Omega))\cap L^2(0,T;H^1(\Omega))}\notag\\
&\quad +\left \|\nabla  \mu_{\epsilon}^{m}\right \|_{L^{2}(0,T;L^2(\Omega))} \le  C_T.
\label{prioria2}
\end{align}
Next, testing \eqref{g4.dr} with $1$, we get
\begin{align}
|\overline{\mu}_{\epsilon}^{m }|&=|\Omega|^{-1}|(\varPsi'_{\epsilon}(\varphi_{\epsilon}^{m }),1)-(\chi \sigma_{\epsilon}^{m },1)| \notag \\
&\le |\Omega|^{-1}\| \varPsi'_{\epsilon}(\varphi_{\epsilon}^{m })\|_{L^1} + C \|  \sigma_{\epsilon}^{m }\|,
\label{average value}
\end{align}
and (cf. \cite{GG,MT})
\begin{equation}
\begin{aligned}
&\| \varPsi'_{\epsilon}(\varphi_{\epsilon}^{m })\|_{L^1}
\le C (1+\|  \nabla \mu_{\epsilon}^{m }(s)\|). 
\end{aligned}
\end{equation}
Using Poincar\'{e}'s inequality and \eqref{intenergy}, we obtain
\begin{equation}
\left \|  \mu_{\epsilon}^{m }  \right \|_{L^{2}(0,T;H^1(\Omega)) }\le C_T.\label{amuL2H1}
\end{equation}
\textbf{Fifth estimate}. Concerning time derivatives, since $ \varphi_{\epsilon}^{m } \in \ L^2(0,T;H^2_{N}(\Omega))$ and $\bm{u}^{m }  \in \  L^{\infty}(0,T;\bm{L}^{\infty}(\Omega))$, similar to \cite{LW} we obtain
\begin{equation}
\left \|  \varphi_{\epsilon}^{m }\bm{u}^{m }\right \|^2_{L^{2}(0,T;\bm{L}^2(\Omega))} \le C,
\label{vphi}
\end{equation}
which together with \eqref{g1.ar} yields
\begin{equation}
\left \|  \partial_{t}\varphi_{\epsilon}^{m }\right \|^2_{L^{2}(0,T;H^1(\Omega)')} \le C.\notag
\end{equation}
Similarly, we also have 
\begin{equation}
\left \|  \partial_{t}\sigma_{\epsilon}^{m }\right \|^2_{L^{2}(0,T;H^1(\Omega)')} \le C.\notag
\end{equation}

\noindent \textbf{Step 3.\ Passage\ to\ the\ limit\ as}\ $\epsilon \to 0$. 
 The estimates obtained in the previous step are independent of $\epsilon$ (nevertheless, they may depend on $m$). Then we are able to pass to the limit as $\epsilon\to 0$ to find a convergent subsequence, using a similar compactness argument like in \cite[Section 4]{MT}. The limit function denoted by $(\varphi^m, \sigma^m)$ is a global weak solution to problem \eqref{g1.a}--\eqref{g6.h}. 
 In particular, it satisfies $\varphi^m\in L^\infty(\Omega\times(0,T))$ such that 
\begin{equation}
-1<\varphi(x,t)<1,\quad  \textrm{a.e. in}\ \Omega\times(0,T). \notag
\ee
Uniqueness of the solution follows from the energy method. We omit the remaining details.

The proof of Lemma \ref{fp} is complete.$\hfill\blacksquare$

\section*{Acknowledgments} 
\noindent The author would like to thank Professors J.-G. Liu and H. Wu for their helpful discussions.


\end{document}